\documentclass[11pt]{amsart}
\usepackage{verbatim}
\usepackage{latexsym}
\usepackage{amsmath,amsthm,amssymb, amscd,color}
\usepackage{mathrsfs}%For mathscr
\usepackage[all]{xy}
\usepackage[all]{xy}
\usepackage{tikz}
\usepackage{tikz-cd}
\usepackage{mathabx}
\usetikzlibrary{decorations.pathmorphing}
\usepackage{anysize}
\marginsize{2.5cm}{2.5cm}{2.5cm}{2.5cm}

\input xy \xyoption{all}
\usepackage{blkarray}  
\usepackage{framed}
\usepackage[utf8]{inputenc}

\numberwithin{equation}{section}
\theoremstyle{plain}
\newtheorem{theorem}{Theorem}[section]
\newtheorem{lemma}[theorem]{Lemma}
\newtheorem{corollary}[theorem]{Corollary}
\newtheorem{proposition}[theorem]{Proposition}

\theoremstyle{definition}
\newtheorem{definition}[theorem]{Definition}
\newtheorem{example}[theorem]{Example}

\newtheorem{remark}[theorem]{Remark}

\newcommand{\QQ}{\mathbb{Q}}
\newcommand{\RR}{\mathbb{R}}
\newcommand{\ZZ}{\mathbb{Z}}

\makeatletter \def\blfootnote{\xdef\@thefnmark{}\@footnotetext}\makeatother

\numberwithin{equation}{section}

\numberwithin{equation}{section}

\begin{document}

\title[Blowdown, $k$-wedge and evenness of quasitoric orbifolds]{Blowdown, $k$-wedge and evenness of quasitoric orbifolds}

%Polytopal $k$-wedge and infinite families of integrally equivariant formal orbifolds

\author[K. Brahma]{Koushik Brahma}
\address{Department of Mathematics, Indian Institute of Technology Madras, India}
\email{koushikbrahma95@gmail.com}

\author[S. Sarkar]{Soumen Sarkar}
\address{Department of Mathematics, Indian Institute of Technology Madras, India}
\email{soumen@iitm.ac.in}

\author[S. Sau]{Subhankar Sau}
\address{Department of Mathematics, Indian Institute of Technology Madras, India}
\email{subhankarsau18@gmail.com}

\subjclass[2010]{57R18, 14M25, 52B11, 55N10}

\keywords{quasitoric orbifold, simplicial toric variety, blowup, blowdown, retraction sequence}

\date{\today}
\dedicatory{}

\abstract 
In this paper, we introduce polytopal $k$-wedge construction and blowdown of a simple polytope and inspect the effect on the retraction sequence of a simple polytope due to $k$-wedge construction and blowdown. In relation to this construction, we introduce the $k$-wedge and blowdown of a quasitoric orbifold. We compare the torsions in the integral cohomologies of $k$-wedges and blowdowns of a quasitoric orbifold with the original one. %Furthermore, we show that the $J$-construction of toric orbifolds in \cite{BSS2} is a particular case of blowdown of quasitoric orbifolds. 
These two constructions provide infinitely many integrally equivariantly formal quasitoric orbifolds from a given one. 

%However, interestingly, the $k$-wedge of a quasitoric orbifold may not be possible to obtain from iterated wedge operation in general.

\endabstract

\maketitle

%\tableofcontents

\section{Introduction}

Simplicial wedge operation is a classical technique in the category of simplicial complexes, see \cite{Ew} and \cite{PB}. The authors of \cite{BBCG10} used this idea in the area of toric topology for the first time. Later, several applications have been exploited in \cite{BBCG15} and \cite{CP2nd}.
 One of the main objectives of these works is to construct infinite families of toric manifolds from a given one which may simplify the presentation of their integral cohomology rings.

Let $K$ be a simplicial complex with the vertex set $\{v_1, \dots, v_m\}$.
%let us fix a vertex $v_i$. 
The simplicial wedge of $K$ on $v_i$ is a simplicial complex with the vertices $\{v_1,\dots, v_{i-1}, v_{i_0}, v_{i_1}, v_{i+1}, ..., v_m\}$ defined by $$K(v_i):=\{v_{i_0}, v_{i_1}\} * \text{link}_K\{v_i\} \cup \{\{v_{i_0}\}, \{v_{i_1}\}\} * (K \setminus \{v_i\})$$
where $*$ implies the join of simplicial complexes.
The dual notion of this construction is called polytopal wedge construction. Precisely,  a simple polytope $P$ is called a polytopal wedge of $Q$ if $K_P$ is a simplicial wedge of $K_Q$, where $K_P, K_Q$ are the dual simplicial complexes of $P, Q$ respectively. We note that $K_P$ is a simplicial complex on the set of codimension-1 faces of $P$. The readers are referred to \cite{BBCG10} and \cite{BP} for details on these concepts. 

On the other hand, Davis and Januszkiewicz introduced toric manifolds and toric orbifolds in the pioneering paper \cite{DJ}. However, they studied several topological properties of toric manifolds. Later, toric orbifolds were explicitly defined in \cite{PS} with the name `quasitoric orbifolds' to avoid similar terminology in algebraic geometry. Weighted projective spaces and simplicial projective toric varieties are some well-known examples of toric orbifolds.  Here, the authors prefer to use the term quasitoric orbifold instead of toric orbifold. A quasitoric orbifold is an even-dimensional effective orbifold equipped with a `locally standard' half-dimensional torus action such that the orbit space has the structure of a simple polytope. The seminal work \cite{Ka} computed the integral cohomology ring of weighted projective spaces. This inspired us to study the integral cohomology of quasitoric orbifolds as this may help in classifying quasitoric orbifolds up to diffeomorphisms. Note that a CW-complex structure can be constructed on an effective orbifold following \cite{Gor}. 
Several works discussed the de-Rham cohomology, the singular cohomology, the Chen-Ruan cohomology ring, orbifold $K$-theory of orbifolds with rational, real or complex coefficients; see \cite[Chapter 2 and 3]{ALR}, \cite{Hat}, \cite{CR}, \cite{ChRu}. However, the computation of these cohomologies with integral coefficients is considerably difficult.

A quasitoric orbifold is called {\it even} if its integral cohomology ring is torsion-free and concentrated in even degrees. The paper \cite{BSS} initiated the investigation of which (quasi)toric orbifold is even. Subsequently,
in \cite{BSS2}, they constructed infinitely many even (quasi)toric orbifolds using the polytopal wedge construction. In this paper, we study several properties of blowdown and $k$-wedge of polytopes and quasitoric orbifolds which generalize the wedge construction of polytopes and $J$-construction of toric orbifolds respectively. Moreover, we extend the discussion on the evenness of quasitoric orbifolds.

The paper is organized as follows. In Section \ref{Sec_Prelim}, we revisit the concept of the retraction sequence (Definition \ref{def:ret_simple_polytope}) of a simple polytope from \cite{BSS}. We recall that a retraction sequence of a simple polytope $Q$ induces retraction sequences of $Q \times \Delta$ for any simplex $\Delta$, see Proposition \ref{retraction of F times Delta}.
Then, following \cite{PS}, we briefly go through the basic construction of a quasitoric orbifold $X(Q,\lambda)$ from a combinatorial data called an $\mathcal{R}$-characteristic pair $(Q, \lambda)$ where $$\lambda \colon \mathcal{F}(Q) \to \ZZ^{\dim Q}$$ is called an $\mathcal{R}$-characteristic function on the simple polytope $Q$, see Definition \ref{toric orbifold construction}. We discuss some invariant subspaces of $X(Q,\lambda)$ corresponding to the faces of $Q$ and the orbifold property of these subspaces. We also recollect the computation of the orbifold singularities at the fixed points of $X(Q, \lambda)$ and its invariant subspaces, see (\ref{singularity}) and (\ref{singularity 1}).

In Section \ref{Sec_Generalized wedge construction of polytopes}, we define polytopal $k$-wedge $Q_F(k)$ of a simple polytope $Q$ at a facet $F$ and prove that $Q_F(k)$ is a simple polytope of dimension $(\dim Q+k)$, see Lemma \ref{Lem_k wedge is a simple polytope}.
We observe that this construction can be carried out at a codimension-$\ell$ face with $2 \leq \ell \leq \dim Q$. However, this may not produce a simple polytope, see Remark \ref{rmk_k_wdg_face}.

In Section \ref{Sec_Blowdown of a simple polytope}, we introduce the concept of blowdown of a convex polytope. We show that the blowdown of a simple polytope may not be a simple polytope in general, see Figure \ref{Fig_blowdown but not simple}.
We also provide the necessary and sufficient conditions when a blowdown preserves the simpleness of a polytope, see Lemma \ref{Q' is simple}.

The main result of this section is that a retraction sequence of $Q$ induces a retraction sequence on its blowdown $Q'$ if $Q'$ is a simple polytope, see Theorem \ref{retraction sequence of blowdown}. Moreover, we construct a retraction sequence of $Q_F(k)$ from a given retraction sequence of $Q$, see Corollary \ref{cor_ret_k_wed} and \ref{cor_k wed_ret}.

In Section \ref{Sec_Blowdown of quasitoric orbifolds},
%we study the blowdowns of quasitoric orbifolds and discuss the evenness of quasitoric orbifolds.
%Let $Q'$ be a blowdown of $Q$ such that $Q'$ is simple. Given two $\mathcal{R}$-characteristic pairs $(Q,\lambda)$ and $(Q',\lambda')$ we define \emph{restriction}, see \eqref{defn of lambda'}. 
first, we define the blowdown of a quasitoric orbifold, see Definition \ref{blowdown definition}.
If $(Q, \lambda)$ and $(Q', \lambda')$ are $\mathcal{R}$-characteristic pairs such that $Q'$ is a blowdown of $Q$ then we analyze when $(Q', \lambda')$ is a restriction of $(Q, \lambda)$ in the sense \eqref{defn of lambda'}. Then, in Theorem \ref{blowdown theorem}, we show that if $(Q,\lambda)$ satisfies some combinatorial conditions along with the hypotheses $(A_2)$ and $(A_3)$ then $(Q',\lambda')$ possesses the similar combinatorial conditions. We show that, in general, we may not be able to remove the hypotheses $(A_2)$ and $(A_3)$ from Theorem \ref{blowdown theorem}; see Example \ref{Eg A2 needed in blowdown} and Example \ref{Eg A3 need in blowdown} respectively.
 We conclude that the integral homology of certain blowdown of a quasitoric orbifold has no $p$-torsion, see 
 Theorem \ref{blowdown no torsion}. If a quasitoric orbifold is obtained by a sequence of blowdown on a quasitoric manifold and each step satiesfies the hypetheses of Theorem \ref{blowdown no torsion} for any prime $p$, then we conclude that the integral cohomology of a blowdown of a quasitoric orbifold is concentrated in even degrees and has no torsion, see Corellary \ref{blowdown free}.

In Section \ref{generalization}, we define $k$-wedge construction on quasitoric orbifolds. We remark that, in general, $k$-wedge on quasitoric orbifold may not be possible to obtain from iterated polytopal wedge construction of \cite{BSS2}.
Also, the blowdown of a simple polytope may not be possible to construct from the polytopal $k$-wedge constructions of a simple polytope, see Example \ref{Eg_blowdown can not be obtained by polytopal wedge construction}.
Consequently, we construct infinitely many integrally equivariantly formal quasitoric orbifolds from a given one in more generality.

\section{Preliminaries}\label{Sec_Prelim}
\subsection{Retraction sequences of polytopes}\label{Subsec_Retraction Sequence of simple polytopes}
In this subsection, we recall a few basics of retraction sequences on polytopes. 
 The convex hull of a finite set of points in $\mathbb{R}^n$ for some $n$ is called a convex polytope.
The vertices, edges, and facets of a convex polytope are faces of dimension $0$, $1$, and $(n-1)$, respectively.
If at each vertex of an $n$-dimensional convex polytope $Q$ exactly $n$ facets intersect, then $Q$ is called a {simple polytope}.
Some well-known examples of simple polytopes are cubes, simplices and prisms.
  We denote the set of vertices of a convex polytope $Q$ by $V(Q)$ and the set of facets of $Q$ by $\mathcal{F}(Q)$ throughout this paper.

\begin{definition}\cite[Definition 5.1]{Zie}
A polytopal complex $\mathcal{C}$ is a finite collection of convex polytopes in $\mathbb{R}^n$ such that the following holds:

\hspace{.8cm} (1) If $E$ is a face of $F$ and $F \in \mathcal{C}$ then $E \in \mathcal{C}$.

\hspace{.8cm} (2) If $E,F \in \mathcal{C}$ and $E \cap F \neq \varnothing$ then $E \cap F$ is a face of both $E$ and $F$.

The dimension of a polytopal complex is defined to be the maximum dimension of the convex polytope in it.
The union of the convex polytopes in $\mathcal{C}$ is called its geometric realization. 
\end{definition}

Let $Q$ be an $n$-dimensional simple polytope and $\mathcal{L}(Q):=\{F \colon F \text{ is a face of } Q\}$. Then $\mathcal{L}(Q)$ is an $n$-dimensional polytopal complex. If $P$ is a subset of $Q$ such that $P$ is the union of some faces of $Q$, then $\mathcal{L}(P)$ is also a polytopal complex.  
For simplicity in this situation, we call $P$ a subcomplex of $Q$.% and denote it by $P \subset Q$. 

\begin{definition}
Let $P$ be a subcomplex of $Q$ and $v \in V(P) \subset V(Q)$. The vertex $v$ is called a free vertex of $P$ if $v$ has a neighborhood $U_v$ in $P$ such that $U_v$ is homeomorphic to $\RR^{d}_{\geq 0}$ as a manifold with corners for some $0 \leq d \leq \dim(P)$. The set $U_v$ is called a local neighborhood of the free vertex $v$ in $P$.
\end{definition}

\begin{definition}\label{def:ret_simple_polytope}
Let $Q$ be a polytope with $m$ vertices and there exists a sequence
$\{(B_{\ell}, E_{\ell}, b_{\ell})\}_{\ell=1}^m $ of triplets such that
\begin{itemize}
\item[(1)] $B_1 = Q=E_1$ and $b_1$ is a free vertex of $Q$.
\item[(2)] $B_{\ell} \subset B_{\ell -1}$ such that $B_{\ell}=\bigcup \{F~ |~F \text{ is a face in } B_{\ell -1} \text{ and } b_{\ell-1} \not \in V(F) \}$.
\item[(3)] $b_\ell$ is a free vertex in $B_\ell$ and $E_{\ell}$ is the maximal dimensional face of $B_{\ell}$ containing the vertex $b_{\ell}$. 
\item [(4)]
$B_m = E_m=b_m$.
\end{itemize}
Then the sequence $\{(B_{\ell}, E_{\ell}, b_{\ell})\}_{\ell = 1}^m$ is called a retraction sequence of $Q$ starting with the vertex $b_1$ and ending at $b_m$. 
\end{definition}

\begin{figure}
\begin{tikzpicture}[scale=0.3]
\draw[dashed] (0,0)--(4,1)--(6,3)--(2,2)--cycle;%behind rectangle
\draw (0,0)--(2,2)--(3,-1)--cycle;%front triangle
\draw[dashed] (4,1)--(6,3)--(7,0)--cycle;
\draw (6,3)--(2,2)--(3,-1)--(7,0)--cycle;
\node at (2,2) {$\bullet$};
\node[above] at (2,2) {$b_1$};

\begin{scope}[xshift=240]
\draw[dashed] (0,0)--(4,1)--(6,3);
\draw (0,0)--(3,-1);
\draw[dashed] (4,1)--(6,3);
\draw[dashed] (4,1)--(7,0);
\draw (6,3)--(7,0);
\draw (3,-1)--(7,0)--cycle;
\node at (6,3) {$\bullet$};
\node[above] at (6,3) {$b_2$};
\end{scope}

\begin{scope}[xshift=480]
\draw[dashed] (0,0)--(4,1);
\draw (0,0)--(3,-1);
\draw[dashed] (4,1)--(7,0);
\draw (3,-1)--(7,0)--cycle;
\node at (4,1) {$\bullet$};
\node[above] at (4,1) {$b_3$};
\end{scope}

\begin{scope}[xshift=720]
\draw (0,0)--(3,-1);
\draw (3,-1)--(7,0)--cycle;
\node at (7,0) {$\bullet$};
\node[above] at (7,-2) {$b_4$};
\end{scope}

\begin{scope}[xshift=960]
\draw (0,0)--(3,-1);
\node at (3,-1) {$\bullet$};
\node[right] at (3,-1) {$b_5$};
\end{scope}

\begin{scope}[xshift=1140]
\node at (0,0) {$\bullet$};
\node[above] at (0,0) {$b_6$};
\end{scope}

\end{tikzpicture}
\caption{A retraction sequence of prism.}
\label{Fig_retraction sequence}
\end{figure}
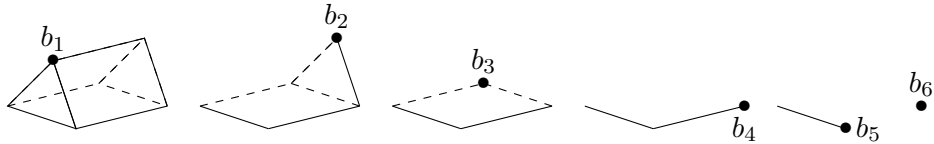

Remark that the conditions (2) and (3) of Definition \ref{def:ret_simple_polytope} imply $B_{\ell} = B_{\ell +1} \cup E_{\ell}$ for $\ell = 1, \ldots, m-1$. Note that a retraction sequence of $Q$ induces an ordering on $V(Q)$. Figure \ref{Fig_retraction sequence} gives an example of a retraction sequence of a prism.
In \cite{BSS}, the authors proved that a simple polytope admits at least one retraction sequence.
We remark that all convex polytopes may not possess retraction sequences in general. But some convex polytopes admit retraction sequences though they are not simple. For example, there is no retraction sequence of the octahedron; however, we can construct a retraction sequence of a pyramid on a pentagonal base.

\begin{proposition}\cite[Proposition 2.5]{BSSau21}\label{retraction of F times Delta}
Let $Q$ be a simple polytope and $\Delta$ be a simplex. Then $Q \times \Delta$ has a retraction sequence induced from the retraction sequences of $Q$ and $\Delta$.
\end{proposition}

\subsection{Some basics of quasitoric orbifolds}
\label{Subsec_Basics of quasitoric orbifold}
A quasitoric orbifold is an even-dimensional effective orbifold with nice enough half-dimensional torus action. We can realize quasitoric orbifolds as a topological analog of simplicial projective toric varieties.  
In this subsection, we briefly recall the constructive definition of a quasitoric orbifold, some notion of invariant suborbifolds, and the singularities at some special points following \cite{PS}.
The authors of \cite{ALR} and \cite{MM} gave a nice exposure to (effective) differentiable orbifolds. 
Let $Q$ be an $n$-dimensional simple polytope with $V(Q):=\{b_1,\dots,b_m\}$ and $\mathcal{F}(Q):=\{F_1,\dots,F_r\}$.

\begin{definition}\label{toric orbifold construction}
Let $\lambda \colon \mathcal{F}(Q) \rightarrow \mathbb{Z}^n$ be a map such that for $i \in \{1, \ldots, r\}$ each $\lambda(F_i)$ is primitive and
\begin{equation}\label{linear independency of char vec}
\{\lambda(F_{i_1}),\dots,\lambda(F_{i_k})\}~ \text{is linearly independent if } \bigcap_{j=1}^k F_{i_j}\neq \varnothing.
\end{equation} 
Then $\lambda$ is called an $\mathcal{R}$-\emph{characteristic function} on $Q$.
The vector $\lambda(F_i)$ is denoted by $\lambda_i$ and called the $\mathcal{R}$-\emph{characteristic vector} assigned to the facet $F_i$.
The pair $(Q,\lambda)$ is called an $\mathcal{R}$-\emph{characteristic pair}.
\end{definition}

\begin{remark}
 Let $F$ be a $d$-dimensional face of an $n$-dimensional simple $Q$ with $d<n$. Then $$F=\bigcap_{j=1}^{n-d} F_{i_j}$$ for some unique facets $F_{i_1},\dots,F_{i_{n-d}}$ of $Q$.
 If the set of vectors $\{\lambda_{i_j}|j=1,\dots,(n-d)\}$ spans an $(n-d)$-dimensional unimodular subspace of $\mathbb{Z}^n$, then  $\lambda$ is called a \emph{characteristic function} and the pair $(Q,\lambda)$ is called a \emph{characteristic pair}, see $(\ast)$ in page 423 of \cite{DJ}.
Note that \cite[Definition 3.5]{SaSu} is a generalization of Definition \ref{toric orbifold construction}.

\end{remark}

\begin{figure}
\begin{tikzpicture}[scale=.4]
\draw (0,0)--(3,0)--(3,3)--(0,3)--cycle;
\node[left] at (3,-.8) {(1,0)};
\node[right] at (3,1.5) {(2,1)};
\node[above] at (1.5,3) {(-3,7)};
\node[left] at (0,1.5) {(5,4)};

\node[left] at (2.5,-3) {(a)};

\begin{scope}[xshift=350]
\draw[dashed] (0,0)--(4,1)--(6,3)--(2,2)--cycle;%behind rectangle
\draw (0,0)--(2,2)--(3,-1)--cycle;%front triangle
\draw[dashed] (4,1)--(6,3)--(7,0)--cycle;
\draw (6,3)--(2,2)--(3,-1)--(7,0)--cycle;

\draw[dotted,thick](5.5,1.5)--(6.33,2);
\draw[->] (6.33,2) to [out=60,in=270] (7,3.5);
\node[above] at (7,3.5) {(1,2,1)};

\draw[dotted,thick] (3,1.5)--(3,2.25);
\draw[->] (3,2.25)--(3,3.5);
\node[above] at (3,3.5) {(1,1,0)};

\draw[->] (1.7,1.2)--(0,2);
\node[left] at (0,2.5) {(1,0,0)};

\draw[dotted, thick] (4.2,0)--(4,-.75);
\draw[->] (4,-.75) to [out=240, in=330] (0,-1.5);
\node[above] at (-1.3,-1.5) {(3,2,1)};

\draw[->] (5,0)to [out=300, in=180](6.5,-1.5); 
\node[right] at (6.5,-1.5) {(0,2,1)};

\node[left] at (5,-3) {(b)};
\end{scope}

\end{tikzpicture}
\caption{Some examples of $\mathcal{R}$-characteristic functions on simple polytopes.}
\label{Fig_characteristic function on simple polytopes}
\end{figure}

\begin{example}
We give an example of an $\mathcal{R}$-{characteristic function} on a square in Figure \ref{Fig_characteristic function on simple polytopes}(a) and on a prism in Figure \ref{Fig_characteristic function on simple polytopes}(b).
\end{example}

We recall the basic construction of a quasitoric orbifold from an $\mathcal{R}$-characteristic pair $(Q, \lambda)$ following \cite{PS}. 
Let $F$ be a face of dimension $d (0 \leq d <n)$ in $Q$. Then $F=\bigcap_{j=1}^{n-d} F_{i_j}$ for some unique facets $F_{i_1}, \dots, F_{i_{n-d}}$ of $Q$. 
Each $\lambda_i \in \ZZ^n$ determines a line in $\mathbb{R}^n(=\ZZ^n \otimes_{\ZZ} \RR)$, whose image under the exponential map $$\exp \colon \mathbb{R}^n \to T^n=(\ZZ^n \otimes_{\ZZ} \RR)/\ZZ^n$$ is a circle subgroup, denoted by $T_i$. Let $T_F:=\left<T_{i_1}, \dots, T_{i_{n-d}}\right>$. Then $T_F$ is an $(n-d)$-dimensional subtorus of $T^n$. We define $T_Q =1 \in T^n$.
Consider the equivalence relation $\sim$ on $T^n \times Q$ is defined by
\begin{equation}\label{eq_equivalence_rel}
(t,x) \sim (s,y)~ \text{if and only if}~ x=y\in \mathring{F} ~\text{and}~ t^{-1}s \in T_F,
\end{equation} 
where $x$ is in the relative interior of the unique face $F$ of $Q$.
The quotient space
\begin{center}
$X(Q,\lambda): =(T^n \times Q) / \sim$.
\end{center}
has an orbifold structure with a natural $T^n$ action. The orbit map
\begin{equation}\label{orbit map}
\pi \colon X(Q,\lambda) \rightarrow Q
\end{equation}
is defined by  $[t,x]_{\sim} \mapsto x$, where $[t,x]_{\sim}$ is the equivalence class of $(t,x)$.
In \cite{PS}, the authors discussed the orbifold structure of the space $X(Q,\lambda)$ explicitly. They also show that the axiomatic definition of quasitoric orbifolds and this definition of quasitoric orbifolds are equivalent. Therefore, studying the topological properties of quasitoric orbifolds with the constructive definition is enough.

Now we discuss the $\mathcal{R}$-{characteristic pairs} for some closed invariant suborbifolds of $X(Q,\lambda)$ following \cite{PS}. Then we compute the singularities of some special points of these invariant suborbifolds.
Consider a $d$-dimensional face $F$ of $Q$ with $0<d < n$. 
Then $F$ is simple and $F= \bigcap_{j=1}^{n-d} F_{i_j}$, for some unique facets $F_{i_1}, \dots, F_{i_{n-d}}$ of $Q$. Let $$N(F):=\left<\lambda_{i_1}, \dots, \lambda_{i_{n-d}}\right>$$ where $\lambda_{i_1}, \dots, \lambda_{i_{n-d}}$ are the $\mathcal{R}$-{characteristic vectors} assigned to these facets respectively. Then $N(F)$ is an $(n-d)$-dimensional submodule of $\mathbb{Z}^n$.

Consider the projection map
\begin{equation}\label{eq_proj_map}
   \rho_F \colon \mathbb{Z}^n \rightarrow \mathbb{Z}^n/{((N(F) \otimes_{\mathbb{Z}} \mathbb{R}) \cap \mathbb{Z}^n)}\cong \ZZ^d. 
\end{equation}

The facets of $F$ are the following
\begin{center}
$ \mathcal{F}(F):=
\{F \cap F_j \mid F_j \in \mathcal{F}(Q) \text{ and } j\neq i_1,\dots,i_{n-d} \text{ and } F \cap F_j \neq \varnothing\}.$
\end{center} 
Then, one can define a map
%$\lambda_F$ on the set of facets of $F$
\begin{align}\label{projection of lambda}
\lambda_F : \mathcal{F}(F) \rightarrow \mathbb{Z}^d
\end{align}
 by $\lambda_F(F \cap F_j):=\text{prim}((\rho_F \circ \lambda)(F_j))$, where $\text{prim}((\rho_F \circ \lambda)(F_j))$ denotes the primitive vector of $(\rho_F \circ \lambda)(F_j)$. 
Note that, $\lambda_F$ is an $\mathcal{R}$-{characteristic function} on $F$. Consequently, it gives a quasitoric orbifold $X(F, \lambda_F)$ which is an invariant suborbifold of $X(Q,\lambda)$, see \cite[Section 2.3]{PS}.

Now we recall how the order of singularities associated to each vertex of the face $F$ is defined. Let $v \in V(F) \subset V(Q)$ and 
\begin{equation}\label{define pi_F}
\pi_{F} \colon X(F, \lambda_F) \to F
\end{equation}
be the orbit map. Then $v=(F \cap F_{j_1}) \cap \dots \cap (F \cap F_{j_d})$ for some unique facets $F_{j_1}, \dots ,F_{j_d}$ of $Q$. The orbifold singularity at the point $\pi_{F}^{-1}(v)$ in $X(F, \lambda_F)$ is defined by
\begin{equation}\label{singularity}
G_{F}(v):= \mathbb{Z}^d / \langle\lambda_F(F \cap F_{j_1}), \dots ,\lambda_F(F \cap F_{j_d})\rangle.
\end{equation}
%Recall that $Q$ is also a face of $Q$.
When $F=Q$, then $v=F_{i_1} \cap \dots \cap F_{i_n}$ for some unique facets $F_{i_j}$ of $Q$. Then the orbifold singularity at the point $\pi^{-1}(v)$ in $X(Q, \lambda)$ is given by
\begin{align}\label{singularity 1}
G_Q(v) := \ZZ^n / \langle\lambda(F_{i_1}), \dots, \lambda(F_{i_n})\rangle.
\end{align}
We call the matrices
\begin{align}\label{Eq_associated matrices}
    A_v^Q&:=\begin{pmatrix}\lambda(F_{i_1})^t & \dots &\lambda(F_{i_n})^t)\end{pmatrix}, ~\mbox{and}\\
    A_v^F&:=\begin{pmatrix}
    \lambda_F(F \cap F_{j_1})^t & \dots & \lambda_F(F \cap F_{j_d})^t\end{pmatrix} \nonumber
\end{align}
associated to the vertex $v$ in $Q$ and $F$ respectively.
Note the following:
\begin{align}\label{singularity 2}
|G_{F}(v)|&=|\det A_v^F|=|\text{det}[\lambda_F(F \cap F_{j_1})^t \dots \lambda_F(F \cap F_{j_d})^t]|, ~ \text{and} \\
|G_Q(v)| &= |\det A_v^Q|=|\text{det}[\lambda(F_{i_1})^t \dots \lambda(F_{i_n})^t]|.\nonumber
\end{align}
The number $|G_F(v)|$ encodes the order of orbifold singularity of the quasitoric orbifold $X(F, \lambda_F)$ at the point $\pi_F^{-1}(v)$.

\section{Polytopal $k$-wedge construction on a simple polytope}\label{Sec_Generalized wedge construction of polytopes}
In this section, we generalize the polytopoal wedge construction in a broader sense and call it $k$-wedge construction on a simple polytope $Q$. We further show that this construction produces another simple polytope of dimension $(\dim Q+k)$.

Let $Q$ be an $n$-dimensionl simple polytope in $\RR^n$ and $F$ a facet of $Q$. 
We consider the polyhedron $Q \times \RR^k_{ \geq 0} \subset  \RR^{n+k}$ and identify $Q \times \textbf{0}_k$ with $Q$ where $\textbf{0}_k$ is the corner $(0, \dots, 0)$ in $\RR_{\geq 0}^k \subset \RR^k$. Let $H$ be a hyperplane in $\RR^{n+k}$ such that it intersects the interior of $Q \times \RR^k_{ \geq 0}$ and divides it into two parts such that one open half space (say $H_{< 0}$) of $H$ contains the vertices $V(Q) \setminus V(F)$ as well as $Q \cap H=F$. 
Let us denote the part containing $Q$ by $Q_F(k)$, that is $$Q_F(k):=(Q \times \RR_{\geq 0}^k) \cap H_{\leq 0}.$$
When $k=1$, the construction is called polytopal wedge construction in \cite{CP} and \cite{CP2nd}.

The hyperplane $H$ can be defined as follows. Choose $n$ many vertices $v_1, \dots, v_n \in V(F)$ which are in `general positions'.
Now we choose $k$-many points $v_{n+1}, \dots, v_{n+k}$ from $v \times \RR^k_{\geq 0}$ such that $v \in V(Q) \setminus V(F)$ and the line segment joining $v \times \textbf{0}_k$ and $v_{n+j}$ is a subset of an edge of $Q \times \RR^k_{\geq 0}$ for $j=1, \dots, k$. Then $\{v_1, \dots, v_n, v_{n+1}, \dots, v_{n+k}\}$ are in general positions in $\RR^{n+k}$. Take the hyperplane $$H:=\big< v_1, \dots, v_n, v_{n+1}, \dots, v_{n+k} \big>.$$
The hyperplane $H$ satisfies the following. Let $p \in F$ be a point and $\overrightarrow{u_p}$ a normal on $F$ towards the interior of $Q$. Then for any $x \in (Q \times \RR^k_{ \geq 0}) \cap H$, the angle between $\overrightarrow{u_p}$ and $\overrightarrow{(x-p)}$ is less than $90^{\circ}$.
Therefore $(v' \times \RR^k_{\geq 0}) \cap H$ is a $(k-1)$-simplex if $ v' \in V(Q) \setminus V(F)$.
%and $(v \times \RR^k_{\geq 0}) \cap H=v$ if $ v \in V(F)$. 
Thus $H$ is a bounding hyperplane for $(Q \times \RR^k_{\geq 0}) \cap H_{\leq 0}$. So $Q_F(k)$ is a convex polytope and
$(Q \times \RR^k_{\geq 0}) \cap H$ is a facet of $Q_F(k)$. Note that $F$ is a face in $Q_F(k)$ of codimension-$(k+1)$, see Figure \ref{Fig_Example of a polytopal $k$-wedge construction.} for an example of this construction.

\begin{figure}
\begin{tikzpicture}[scale=0.7]
    \draw (0,0)--(3,0)--(4,1)--(1,1)--cycle;
    \draw[thick, blue] (0,0)--(1,1);
    \node at (0,0) {$\bullet$};
    \node at (1,1) {$\bullet$};
    \node[left] at (0,0) {$v_1$};
    \node[above] at (1,1) {$v_2$};
    \node at (2,-.75) {$Q$};
     \node[right] at (3,-.2) {$v$};
    \node[left] at (.5,.75) {$F$}; 

    \begin{scope}[xshift=200]
    \draw (0,0)--(3,0)--(4,1)--(1,1)--cycle;
    \draw[thick, blue] (0,0)--(1,1);
    \draw (0,0)--(0,3);
    \draw (3,0)--(3,3);
    \draw (4,1)--(4,4);
    \draw (1,1)--(1,4);
    \draw [fill=yellow, opacity=.4] (0,0)--(3,2)--(4,3.5)--(1,1)--cycle;
    \node at (0,0) {$\bullet$};
    \node at (1,1) {$\bullet$};
    \node at (3,2) {$\bullet$};
    \node[left] at (0,0) {$v_1$};
    \node[above] at (1.3,1) {$v_2$};
    \node[right] at (3,2) {$v_3$};
      \node[right] at (3,-.2) {$v$};
    \node at (2,-.75) {$Q \times \RR_{\geq 0}$};
    \node[left] at (.75,.8) {$F$}; 
    \draw[->] (2,4)--(2.5,2.5); 
    \node[above] at (2,4) {$Q \times \RR_{\geq 0} \cap H$};
    \end{scope}

    \begin{scope}[xshift=400]
    \draw (0,0)--(3,0)--(4,1);
    \draw[thick, dashed] (1,1)--(4,1);
    \draw[thick, blue] (0,0)--(1,1);
    \draw (3,0)--(3,2);
    \draw (4,1)--(4,3.5);
    \draw (0,0)--(3,2)--(4,3.5)--(1,1)--cycle;
    \node at (0,0) {$\bullet$};
    \node at (1,1) {$\bullet$};
    \node at (3,2) {$\bullet$};
    \node[left] at (0,0) {$v_1$};
    \node[left] at (1,1) {$v_2$};
    \node[right] at (3,2) {$v_3$};
      \node[right] at (3,-.2) {$v$};
    \node at (2,-.75) {$Q_F(1)$};
    \end{scope}
\end{tikzpicture}
\caption{Example of a polytopal $k$-wedge construction.}
\label{Fig_Example of a polytopal $k$-wedge construction.}
\end{figure}

\begin{lemma}\label{Lem_k wedge is a simple polytope}
Let $Q$ be an $n$-dimensional simple polytope with a facet $F$. Then  $Q_F(k)$ is an $(n+k)$ dimensional simple polytope. 
\end{lemma}
\begin{proof}
By definition, $Q_F(k)$ is a convex polytope. Thus it is enough to show that, at every vertex, exactly $n+k$ facets of $Q_F(k)$ intersect. 
Note that, the polyhedron $$L^{k-1}_s := \{(x_1, \dots, x_k)\in \RR^k_{\geq 0} | x_s=0\}$$ is a facet of $\RR^k_{\geq 0}$ for $s=1, \dots, k$.
Let 
\begin{align*}
V(Q)=\{v^Q_1, \dots, v^Q_m\} \quad \text{and} \quad
\mathcal{F}(Q)=\{F^Q_1, \dots, F^Q_r\}.
\end{align*}
be the vertex and facet set of $Q$, respectively.
Then the facets of $Q \times \RR^k_{\geq 0}$ are given by 
\begin{align*}
\{Q \times L^{k-1}_s ~| ~ s=1, \dots, k\} \cup
\{F^Q_j \times \RR^k_{\geq 0}~ | ~ j=1, \dots, r\}.
\end{align*} 
Without loss of generality, let $F^Q_r=F$. Then the facet set of $Q_F(k)$ is given by
\begin{equation}\label{Eq_facet set of generalized polytopal wedge}
\mathcal{F}(Q_F(k))=\{F_1, \dots, F_r, F_{r+1}, \dots, F_{r+k}\}
\end{equation}
where 
\begin{equation*} 
F_i:=
\begin{cases}
(F^Q_i \times \RR^k_{\geq 0}) \cap H_{\leq 0}  \quad & \text{ for } i=1, \dots, r-1\\
(Q \times \RR^k_{\geq 0}) \cap H \quad & \text{ for } i=r \\
(Q \times L^{k-1}_s)  \cap H_{\leq 0} \quad & \text{ for } i=r+s, ~s=1,\dots,k. 
\end{cases}
\end{equation*}
Note that there exists a projection $\rho \colon F_r= (Q \times \RR^k) \cap H \to Q $. So $\rho$ is face preserving, and it takes facets to facets. Also, naturally, $Q$ is identified with the face $Q \times {\bf 0}_k$ of $Q_F(k)$.

Let $\textbf{v} \in V(Q_F(k)) \setminus V(Q)$. 
Then $\rho(\textbf{v})=v^Q$ for some $v^Q=\bigcap_{j=1}^n F^Q_{i_j} \in V(Q)$ where $F^Q_{i_j}$'s are some unique facets in $\mathcal{F}(Q) \setminus \{F^Q_r\}.$ 
Then we have $$\textbf{v}=(\bigcap_{j=1}^n F_{i_j}) \cap F_r \cap (\bigcap_{\substack{s=1\\s \neq t}}^k F_{r+s})$$ for some $t \in \{1, ..., k\}$. Thus, exactly $(n+k)$ facets intersect at $\textbf{v}$ in $Q_F(k)$ for this case.

Let $\textbf{v} \in V(Q)\setminus V(F)$. Considering $\textbf{v}$ as a vertex of the simple polytope $Q$ and denote it by $v^Q$, we have $v^Q=\bigcap_{j=1}^n F^Q_{i_j}$ for some unique facets $F^Q_{i_1}, ..., F^Q_{i_n}$ in $\mathcal{F}(Q) \setminus \{F^Q_r\}.$ Thus,  for this case, we have $$\textbf{v}=(\bigcap_{j=1}^n F_{i_j}) \cap (\bigcap_{s=1}^k F_{r+s}).$$ Therefore, $\textbf{v}$ is the intersection of $(n+k)$ many facets in $Q_F(k)$.

Let $\textbf{v} \in V(F)$. Then $\textbf{v}=v^Q \in V(F) \subset V(Q) \subset V(Q_F(k))$. So $v^Q=(\bigcap_{j=1}^{n-1} F^Q_{i_j}) \cap F^Q_r$ for some unique facets $F^Q_{i_1}, ..., F^Q_{i_{n-1}}$ of $Q$. Thus,
considering ${\bf v}$ as a vertex of $Q_F(k)$, we get
$$\textbf{v}=(\bigcap_{j=1}^{n-1} F_{i_j}) \cap F_r \cap (\bigcap_{s=1}^k F_{r+s}).$$ Therefore, in this case also, $\textbf{v}$ is the intersection of $(n+k)$ many facets of $Q_F(k)$.
Thus we get the result.
\end{proof}

We call the simple polytope $Q_F(k)$ the \emph{polytopal $k$-wedge} of $Q$ at $F$.

\begin{figure}
\begin{tikzpicture}[scale=0.7]
    \draw (0,0)--(3,0)--(4,1)--(1,1)--cycle;
    \node at (0,0) {$\bullet$};
    \node[left] at (0,0) {$v_1$};
    \node at (2,-.75) {$Q$};

    \begin{scope}[xshift=200]
    \draw (0,0)--(3,0)--(4,1)--(1,1)--cycle;
    \draw (0,0)--(0,3);
    \draw (3,0)--(3,3);
    \draw (4,1)--(4,4);
    \draw (1,1)--(1,4);
    \draw [fill=yellow, opacity=.4] (0,0)--(3,2)--(4,3.75)--(1,2)--cycle;
    \node at (0,0) {$\bullet$};
    \node at (1,2) {$\bullet$};
    \node at (3,2) {$\bullet$};
    \node[left] at (0,0) {$v_1$};
    \node[left] at (1,2) {$v_2$};
    \node[right] at (3,2) {$v_3$};
    \node at (2,-.75) {$Q \times \RR_{\geq 0}$};
    \node[left] at (.75,.8) {$F$}; 
    \draw[->] (2,4)--(2.5,2.5); 
    \node[above] at (2,4) {$Q \times \RR_{\geq 0} \cap H$};
    \end{scope}

    \begin{scope}[xshift=400]
    \draw (0,0)--(3,0)--(4,1);
    \draw[thick, dashed] (1,1)--(4,1);
    \draw (0,0)--(1,2);
    \draw[thick, dashed] (0,0)--(1,1);
    \draw[thick, dashed] (1,2)--(1,1);
    \draw (3,0)--(3,2);
    \draw (4,1)--(4,3.75);
    \draw (0,0)--(3,2)--(4,3.75)--(1,2)--cycle;
    \node at (0,0) {$\bullet$};
    \node at (1,2) {$\bullet$};
    \node at (3,2) {$\bullet$};
    \node[left] at (0,0) {$v_1$};
    \node[left] at (1,2) {$v_2$};
    \node[right] at (3,2) {$v_3$};
    \node at (2,-.75) {$Q_{v_1}(1)$};
    \end{scope}
\end{tikzpicture}
\caption{}
\label{Fig_k-wedge not simple polytope}
\end{figure}

\begin{example}
Let $Q$ be an interval $I=[0,1]$ with two facets $\{0\}$ and $\{1\}$. If we take $k=2$ and $F=\{1\}$ then the polytopal $2$-wedge of $Q$ at $F$ is a tetrahedron. In Figure \ref{Fig_Example of a polytopal $k$-wedge construction.}, we provide another example.
\end{example}

\begin{remark}\label{rmk_k_wdg_face}
If $F$ is not a facet in $Q$, then the construction $Q_{F}(k)$ may not give a simple polytope in general.
Consider a square $Q$ and a vertex $v_1 \in V(Q)$ as in Figure \ref{Fig_k-wedge not simple polytope}. Then $Q \times \RR_{\geq 0} \subset \RR^3$. Next we take $v_2$ and $v_3$ in general positions of $Q \times \RR_{\geq 0}$ to construct the hyperplane $H$ and eventually $Q_{v_1}(1)$. Note that at $v_1 \in Q_{v_1}(1)$, four facets intersect while $Q_{v_1}(1)$ is $3$-dimensional, see Figure \ref{Fig_k-wedge not simple polytope}. Thus $Q_{v_1}(1)$ is not a simple polytope.
\end{remark}

\section{Blowdowns of polytopes}
\label{Sec_Blowdown of a simple polytope} 
The concept of blowdown of a simple polytope was discussed in \cite[Section 4]{GP2} as follows.
If $Q_1$, $Q_2$ are simple polytopes and $Q_1$ is a blowup (see Definition \ref{defn of blowup of simple polytope}) of $Q_2$ then $Q_2$ is called a blowdown of $Q_1$. But it is not a precise definition, as $Q_1$ may be a blowup of another simple polytope $Q_3$, see Figure \ref{Fig_blowdown as in gp}. In this section, we give the precise definition of blowdown of convex polytope which enriches its beauty.
\begin{figure}
\begin{tikzpicture}[scale=.4]
\draw (2,4)--(1,5)--(-1,5)--(-2,4)--cycle;
\draw (-2,4)--(-2,1)--(2,1)--(2,4);
\draw [dashed] (-2,1)--(-1,2)--(1,2)--(2,1);
\draw [dashed] (-1,2)--(-1,5);
\draw [dashed] (1,2)--(1,5);
\draw [very thick] (2,1)--(2,4);
\draw [<-] (2.1,3.4)--(3,4.5);
\node at (0,-1) {$Q_2$};
\node [right] at (2.8,4.7) {$F_2$};

\begin{scope}[xshift=250]
\draw (0,3)--(2,4)--(1,5)--(-1,5)--(-2,4)--cycle;
\draw (-2,4)--(-2,1)--(0,0)--(2,1)--(2,4);
\draw (0,0)--(0,3);
\draw [dashed] (-2,1)--(-1,2)--(1,2)--(2,1);
\draw [dashed] (-1,2)--(-1,5);
\draw [dashed] (1,2)--(1,5);
\draw [fill=yellow, opacity=.4] (0,0)--(2,1)--(2,4)--(0,3)--cycle;
\node at (0,-1) {$Q_1$};
\draw[<-] (-5.5,2.5)--(-2.5,2.5);
\draw[->] (2.5,2.5)--(5.5,2.5); 
\end{scope}

\begin{scope}[xshift=500]
\draw (1,5)--(-1,5)--(-2,4);
\draw (-2,4)--(-2,1)--(0,0)--(2,1);
\draw [dashed] (-2,1)--(-1,2)--(1,2)--(2,1);
\draw [dashed] (-1,2)--(-1,5);
\draw [dashed] (1,2)--(1,5);
\draw (-2,4)--(0,0);
\draw (1,5)--(2,1);
\draw [very thick] (0,0)--(2,1);
\draw [<-] (1.2,.5)--(2.5,.5);
\node [right] at (2.5,.4) {$F_3$};
\node at (0,-1) {$Q_3$};
\end{scope}

\end{tikzpicture}
\caption{Both $Q_2$ and $Q_3$ are blowdowns as in \cite{GP2}.}\label{Fig_blowdown as in gp}
\end{figure}
We provide the necessary and sufficient conditions for which a blowdown of a simple polytope is again a simple polytope. We also show that the new one possesses an `induced retraction sequence' in the sense of Definition \ref{def:ret_simple_polytope}. 
We prove that blowdown of a polytope is a generalization of the polytopal $k$-wedge construction which is a generalization of the polytopal wedge construction of \cite{CP}.

\begin{definition}[Blowup of a convex polytope]\label{defn of blowup of simple polytope}
Let $Q$ be an $n$-dimensional convex polytope in $\mathbb{R}^n$ and $F$ be a face of $Q$. Take an $(n-1)$ dimensional hyperplane $H$ in $\mathbb{R}^n$ such that one open half space (say $H_{<0}$) contains $V(F)$ and $V(Q) \setminus V(F)$ is a subset of the other open half space $H_{>0}$. Then $\widebar{Q}:=Q \cap H_{\geq 0}$ is called a \emph{blowup} of $Q$ along the face $F$.
\end{definition}

Note that if $F_i$ is a facet of $Q$, then $\widetilde{F}_i:=F_i \cap \widebar{Q}$ is the facet of $\widebar{Q}$ corresponding to $F_i$.
The new facet $\widebar{Q} \cap H$ is called the facet corresponding to the face $F$ and denoted by $\widetilde{F}$.
We refer the reader to \cite{Joy} for several properties of manifold with corners and maps between them.

\begin{definition}[Blowdown of a convex polytope]\label{defn of blowdown of simple polytope}
 
 Let F and $\widetilde{F}$ be two faces of an $n$-dimensional convex polytope $Q$ such that $F \subset \widetilde{F}$ and $\widetilde{F}$ is a facet. Let $Q'$ be a convex polytope with a face $F'$ such that $F'$ is homeomorphic to $F$ as a manifold with corners.
 If the blowup $\widebar{Q'}$ of $Q'$ along the face $F'$ is homeomorphic to $Q$ as a manifold with corners and the restriction on the facet $\widetilde{F'}$ is homeomorphic to $\widetilde{F}$ as a manifold with corners where $\widetilde{F'}$ is corresponding to the face $F'$, then $Q'$ is called a \emph{blowdown} of $Q$ of the facet $\widetilde{F}$ on $F$.
\end{definition}

\begin{figure}
\begin{tikzpicture}[scale=.4]

\draw (0,3)--(2,4)--(1,5)--(-1,5)--(-2,4)--cycle;
\draw (-2,4)--(-2,1)--(0,0)--(2,1)--(2,4);
\draw (0,0)--(0,3);
\draw [dashed] (-2,1)--(-1,2)--(1,2)--(2,1);
\draw [dashed] (-1,2)--(-1,5);
\draw [dashed] (1,2)--(1,5);
\draw [fill=yellow, opacity=.4] (0,0)--(2,1)--(2,4)--(0,3)--cycle;
\node at (0,-1) {$Q$};
\draw [<-] (1.2,.5)--(2.5,.5);
\node [right] at (2.5,.4) {$F$};
\node at (.5,1.3) {$\widetilde{F}$};
%\draw[<-] (-5.5,2.5)--(-2.5,2.5);
\draw[->] (2.5,2.5)--(5.5,2.5);

\begin{scope}[xshift=250]
\draw (1,5)--(-1,5)--(-2,4);
\draw (-2,4)--(-2,1)--(0,0)--(2,1);
\draw [dashed] (-2,1)--(-1,2)--(1,2)--(2,1);
\draw [dashed] (-1,2)--(-1,5);
\draw [dashed] (1,2)--(1,5);
\draw (-2,4)--(0,0);
\draw (1,5)--(2,1);
%\draw [very thick] (0,0)--(2,1);
\draw [<-] (1.2,.5)--(2.5,.5);
\node [right] at (2.5,.4) {$F$};
\node at (0,-1) {$Q'$};
\draw[->] (2.5,2.5)--(5.5,2.5); 
\end{scope}

\begin{scope}[xshift=500]
\draw (1,5)--(-1,5)--(-2,4);
\draw (-2,4)--(-2,1)--(0,0)--(2,1);
\draw [dashed] (-2,1)--(-1,2)--(1,2)--(2,1);
\draw [dashed] (-1,2)--(-1,5);
\draw [dashed] (1,2)--(1,5);
\draw (-2,4)--(0,2);
\draw (1,5)--(2,3);
%\draw [very thick] (0,0)--(2,1);
\draw [<-] (1.2,.5)--(2.5,.5);
\draw [fill=yellow, opacity=.4] (0,0)--(2,1)--(2,3)--(0,2)--cycle;
\node [right] at (2.5,.4) {$F$};
\node at (.5,1.3) {$\widetilde{F'}$};
\node at (0,-1) {$\widebar{Q'}$};
\end{scope}
\end{tikzpicture}
\caption{Blowdown of a pentagonal prism of the face $\widetilde{F}$ on $F$.}\label{Fig_Blowdown}
\end{figure}

\begin{remark}\label{define f}
If $f\colon Q \to \widebar{Q'}$ is the homeomorphism as a manifold with corners in Definition \ref{defn of blowdown of simple polytope}, then $f|_{\widetilde{F}} \colon \widetilde{F} \to \widetilde{F'}$ is a homeomorphism as a manifold with corners. 
Also let $\mathcal{F}(Q)=\{\widetilde{F}, F_1, \dots, F_r\}$ be the facets of $Q$ and $F$ a proper face of $\widetilde{F}$.
Then the facets $F'_1, \ldots, F'_r$ of $Q'$ are such that the facet $\widetilde{F'_i}$ in $\widebar{Q'}$ is homeomorphic to $F_i$ as a manifold with corners through $f|_{F_i}$ for $1 \leq i \leq r$.
\end{remark}

\begin{example}
Let $Q$ be a pentagonal prism with $\widetilde{F}$ and $F$ as shown in Figure \ref{Fig_Blowdown}(a). Also, $Q'$ and $F'$ be as in Figure \ref{Fig_Blowdown}(b) where $F \cong F'$ is a manifold with corners. The blowup of $Q'$ along $F'$ is $\widebar{Q'}$ in Figure \ref{Fig_Blowdown}(c), which is homeomorphic to $Q$ as a manifold with corners. So $Q'$ is a blowdown of $Q$ of the facet $\widetilde{F}$ on $F$.
\end{example}

We note that a blowdown of a simple polytope may not be simple in general, see Figure \ref{Fig_blowdown but not simple}. However, the following lemma gives a criterion when a blowdown of a simple polytope is simple.

\begin{lemma}\label{Q' is simple}
Let $Q$ be an $n$-dimensional simple polytope having a facet $\widetilde{F}$ homeomorphic to $F \times \Delta^{n-d-1}$ as a manifold with corners where $F$ is a face of $\widetilde{F}$ and $\Delta^{n-d-1}$ is a simplex for $0\leq \dim(F)= d \leq (n-1)$. Let $Q'$ be the blowdown of $Q$ of the facet $\widetilde{F}$ on $F$. Then $Q'$ is an $n$-dimensional simple polytope. 
Conversely if $Q'$, the blowdown of $Q$ of the facet $\widetilde{F}$ on $F$, is simple then $\widetilde{F}$ is homeomorphic to $F \times \Delta^{n-d-1}$ as a manifold with corners.
\end{lemma}

\begin{proof}
Let $V(F)=\{b_{\ell_1},b_{\ell_2}, \dots,b_{\ell_k}\}$ and $V(\Delta^{n-d-1})=\{v_1, \dots, v_{n-d}\}$.
Then we may write
\begin{equation}
V(\widetilde{F}):=\{(b_{\ell_i},v_q)\colon 1\leq i \leq k \text{ and }1\leq q \leq (n-d) \}
 \subset V(Q).
\end{equation}
By the definition of blowdown, $Q$ is homeomorphic to $\widebar{Q'}$ as a manifold with corners. We denote this homeomorphism by $f \colon Q \to \widebar{Q'}$ as in Remark \ref{define f}. This induces a bijection between $V(Q)$ and $V(\widebar{Q'})$.

\begin{figure}
\begin{tikzpicture}[scale=.4]
\draw (0,0)--(3,0)--(3,3)--(0,3)--cycle;
\draw (0,0)--(0,3)--(-1.5,4.5)--(-1.5,1.5)--cycle;
\draw (-1.5,4.5)--(1.5,4.5)--(3,3);
\draw [dashed] (-1.5,1.5)--(1.5,1.5)--(3,0);
\draw [dashed] (1.5,1.5)--(1.5,4.5);

\node at(3,3) {$\bullet$};
\draw[->] (-4,4) to [out=45,in=120] (1,3.5);
\draw[->] (4,3.5)--(3.2,3.2);
\node [right] at (4,3.5) {$F$};
\node [left] at (-3.8,3.5) {$\widetilde{F}$};
\node [right] at (1,-.8) {$Q$};
\draw [->] (3.5,2)--(11.5,2);

\begin{scope}[xshift=400]
\draw (0,0)--(3,0)--(3,3)--cycle;
\draw (0,0)--(-1.5,1.5)--(3,3)--cycle;
\draw[dashed] (-1.5,1.5)--(1.8,1.4)--(3,0);
\draw[dashed] (1.8,1.4)--(3,3);

\node [right] at (1,-.8) {$Q'$};
\node at(3,3) {$\bullet$};
\node [right] at (3,3.5) {$F'$};
\end{scope}

\end{tikzpicture}
\caption{Blowdown of a simple polytope may not be a simple polytope.}
\label{Fig_blowdown but not simple}
\end{figure}

Recall that $Q' \cap H_{\geq 0} =\widebar{Q'}$, see Definition \ref{defn of blowup of simple polytope}. Then $M_g=Q' \cap H_{\geq 0}$ is a mapping cylinder for the projection map \[g\colon \widetilde{F'} \xrightarrow{\cong} F' \times \Delta^{n-d-1} \to F'.\] So there is a face preserving homotopy of $M_g$ on $F'$. 
Let us consider a tubular neighborhood $N_{\widetilde{F'}}$ of $\widetilde{F'}$ in $\widebar{Q'}$ such that it does not contain any vertices in $V(\widebar{Q'}) \setminus V(\widetilde{F'})$. We define 
\begin{equation}
f' \colon \widebar{Q'} \to Q'
\end{equation}
 by $f'(N_{\widetilde{F'}}) \longmapsto N_{\widetilde{F'}} \cup M_g$ preserving the face structure and $f'(x)=x$ if $x \in \widebar{Q'} \setminus N_{\widetilde{F'}}$. Now let 
\begin{equation}\label{define tilde f}
\tilde{f}=f' \circ f \colon Q \to Q'.
\end{equation}
Since $f$ is a homeomorphism and face preserving, the map $\tilde{f}$ is a face preserving map. Then
\begin{equation}
V(Q')=\{\tilde{f}(v):v \in V(Q) \}
\end{equation}
where $\tilde{f}(b_{\ell_i},v_1)=\tilde{f}(b_{\ell_i},v_2)= \dots = \tilde{f}(b_{\ell_i},v_{n-d})$ for all $i=1, \dots ,k$,  and $\tilde{f}$ is one-one otherwise.

Let $\tilde{f}(b) \in V(Q')$ such that $b \in V(\widetilde{F})$. Then $b=(b_{\ell_i},v_q)$ for some $b_{\ell_i} \in V(F)$ and ${v_q \in V(\Delta^{n-d-1})}$ and $\widetilde{f}(b)=\widetilde{f}(b_{\ell_i})$.
Then, considering $b_{\ell_i}$ as a vertex of $Q$, we have 
\begin{equation}\label{b_l_i}
b_{\ell_i} = (\bigcap_{j=1}^{n-1} F_{i_j}) \cap \widetilde{F}
\end{equation}
for some unique facets $F_{i_1}, ..., F_{i_{n-1}}$ of $Q$.
 Also $\{b_{\ell_i}\} \times \Delta^{n-d-1}$ is homeomorphic to $\Delta^{n-d-1}$ as a manifold with corners.
Let $\mathcal{F}(\Delta^{n-d-1}):=\{F_1^\Delta, \dots, F_{n-d}^\Delta\}$ be the set of facets of $ \Delta^{n-d-1}$. 
Exactly one of these facets does not contain $b_{\ell_i}$, say $F_1^{\Delta}$ without loss of generality.
If we consider $F \times F_s^\Delta$ for $1 \leq s \leq (n-d)$ then they are facets of $\widetilde{F}$ and codimension $2$ faces in $Q$.
Note that these may not be the total collection of facets of $\widetilde{F}$.
Thus $F \times F_s^\Delta = \widetilde{F} \cap P_{i_s}$ for a unique facet $P_{i_s}$ in $Q$ for all $1 \leq s \leq (n-d)$.
Except $P_{i_1}$ all other facets in $\{P_{i_s} \colon 2 \leq s \leq (n-d)\}$ contain $b_{\ell_i}$. So 
\begin{equation}\label{subset for blowdown char eqn}
\{P_{i_s} \colon 2 \leq s \leq (n-d)\} \subseteq \{F_{i_j} \colon 1 \leq j \leq n-1\}.
\end{equation}
Since $\widetilde{f}$ is face preserving and $\tilde{f}(\widetilde{F})=F$, from Remark \ref{define f} and \eqref{b_l_i}, we have
\begin{equation}\label{tilde f b_l_i}
\widetilde{f}(b) =(\bigcap_{j=1}^{n-1} \widetilde{f}(F_{i_j})) \cap F =(\bigcap_{j=1}^{n-1} \widetilde{f}(F_{i_j})) \cap \widetilde{f}(P_{i_1})
\end{equation}
where $P_{i_1}$ is described in the previous paragraph. So at the vertex $\widetilde{f}(b)$ in $Q'$ exactly $n$ facets intersect.
A similar construction can be done for any vertex of $F$.

Now let $b \in Q \setminus \widetilde{F}$ be any vertex. Then $b = \bigcap_{t=1}^n F_{i_t}$ for some unique facets $F_{i_t}$ of $Q$ and 
\begin{equation}\label{vertex not in facet}
\widetilde{f}(b) = \bigcap_{t=1}^n \widetilde{f}(F_{i_t}).
\end{equation}
This concludes at every vertex of $Q'$ exactly $n$ facets meet. So, $Q'$ is a simple polytope.

The converse part follows from Definition \ref{defn of blowup of simple polytope} and \ref{defn of blowdown of simple polytope} as $Q$ and $\widebar{Q'}$ are homeomorphic as manifold with corners. 
Precisely, since $Q'$ is simple and $F(\cong F' \text{ as manifold with corners})$ is a face then the facet $\widetilde{F'}$ of $\widebar{Q'}$ corresponding to face $F'$ is homeomorphic to $F \times \Delta$ as a manifold with corners where $\Delta$ is a simplex of dimension $(\dim(Q) -\dim(F) -1)$.
\end{proof}

\begin{remark}\label{blowdown of face}
Let $Q'$ be a blowdown of $Q$ of the facet $\widetilde{F}$ on $F$ such that $Q'$ is simple. 
 If $E \cap \widetilde{F} = \varnothing$, for a face $E$ of $Q$ then $\widetilde{f}(E)$ is homeomorphic to $E$ as manifold with corners.  
\end{remark}

The following lemma investigates how a face $E$ of $Q$ is changed due to blowdown when $E \cap \widetilde{F} \neq \varnothing$. Note that $E \cap \widetilde{F}$ is a face of $\widetilde{F} \cong F \times \Delta^{n-d-1}$.
Thus $E \cap \widetilde{F} = E^F \times E^{\Delta}$ as manifold with corners for some faces $E^F$ and $E^{\Delta}$ of $F$ and $\Delta^{n-d-1}$ respectively.
Now $E^{\Delta}$ is again a simplex as it is a face of the simplex $\Delta^{n-d-1}$. Thus $$E \cap \widetilde{F} = E^F \times \Delta^{q}$$ for some $0 \leq q \leq (n-d-1)$. The face $E \cap \widetilde{F}$ shrinks to $E^F$ due to blowdown.

\begin{lemma}\label{lem:blowdown of face}
 Let $Q'$ be a blowdown of $Q$ of the facet $\widetilde{F}$ on $F$ such that $Q'$ is simple. 
 If $E \cap \widetilde{F} \neq \varnothing$, then $\widetilde{f}(E)$ is either a blowdown of $E$ of the facet $E \cap \widetilde{F}$ on $E^F$ or homeomorphic to a face of $E$ as a manifold with corners.
\end{lemma}
\begin{proof}
Let $\dim(E \cap \widetilde{F})=\beta$. If $\beta=0$, i.e., $E \cap \widetilde{F}$ is a vertex, then $\widetilde{f}(E)$ is homeomorphic to $E$ as manifold with corners. 
Now we consider the cases where $0 < \beta \leq (n-1)$. It is evident from $\dim(E \cap \widetilde{F})=\beta$ that $\dim(E) \geq \beta$. If $\dim(E)=\beta$, then $E \subset \widetilde{F}$ and $\widetilde{f}(E)$ is homeomorphic to a face of $E$.
If $\dim(E)=\beta+1$, then $E \cap \widetilde{F}$ is a facet of $E$. We have $E \cap \widetilde{F} = E^F \times {\Delta}^q$. If $q=0$, then $\widetilde{f}(E)$ is homeomorphic to $E$ as a manifold with corners. If $q>0$ then $\widetilde{f}(E)$ is a blowdown of $E$ of the face $E \cap \widetilde{F}$ on $E^F$.

Now we show $\dim(E) \neq \beta +j$ for $j \geq 2$. First, let $\dim(E)=\beta +2$ and $v \in E \cap \widetilde{F}$ be a vertex. 
Exactly $\beta+2$ edges meet at $v$ in $E$, out of which $\beta$ edges are also edges of $\widetilde{F}$ and two are not. Also, $v$ being a vertex of the facet $\widetilde{F}$, exactly $n-1$ edges meet at $v$ in $\widetilde{F}$. 
This implies that exactly $n+1$ edges meet at $v$ in $Q$, which is a contradiction to $Q$ is an $n$-dimensional simple polytope. Therefore, $\dim(E) \neq \beta+2$, and by similar observation, $\dim(E) \neq \beta+j$ for $j > 2$. 
Thus, the claim of the lemma follows. 
\end{proof}

%The Lemma \ref{Q' is simple} states equivalency of the following statements:
%\begin{enumerate}
%\item The blowdown of $Q$, $Q'$ is simple.
%\item $Q'$ is the blowdown of $Q$ of the face $\widetilde{F}$ on $F$ where $\widetilde{F} \cong F \times \Delta^{n-d-1}$ as a manifold with corners.
%\end{enumerate}
%Observe that the proof of Lemma \ref{Q' is simple} can also be done by counting the edges at each vertex and proving the edges are linearly independent at each vertex.
%We are going to use the concept of $\widetilde{f}$ in \eqref{define tilde f} in our future proofs.
%Note that Remark \ref{remark after blowup}(1) says that the converse of Lemma \ref{Q' is simple} is true. That is if $Q'$ is the blowdown of $Q$ of the face $\widetilde{F}$ on $F$ and $Q'$ is simple then $\widetilde{F}$ is homeomorphic to $F \times \Delta$ as a manifold with corners where $\Delta$ is a simplex of dimension $(\dim(Q) -\dim(F) -1)$.

%Let $P$ be an $(n-1)$-dimensional simple polytope and  $F$  a facet of $P$. Then there exists a blowdown $(P \times [0,1])'$ of the simple polytope $P \times [0, 1]$ of the facet $F \times [0, 1]$ on $F$.

\begin{corollary}\label{cor_gen of poly wedge cons}
Polytopal $k$-wedge of a simple polytope $Q$ at a facet $F$ is a blowdown of $Q \times \Delta^k$ of the facet $F \times \Delta^k$ on $F$.
\end{corollary}

\begin{proof}
If we blowup $Q_F(k)$ along the face $F$, then $\widebar{Q_F(k)}$ is homeomorphic to $Q\times \Delta^k$ as manifold with corners where $\Delta^k$ is a $k$-simplex.
Also the facet $F \times \Delta^k$ arises in $\widebar{Q_F(k)}$ corresponding to the facet $F$ of $Q$. 
Thus the polytopal $k$-wedge construction of $Q$ at $F$ is nothing but a blowdown of $Q \times \Delta^k$ of the face $F \times \Delta^k$ on $F$.
\end{proof}

%There are some blowdowns where $\widetilde{F}$ need not be homeomorphic to $F \times I$ as a manifold with corners to conclude $Q'$ is a simple polytope.
%  For example, if $Q$ is a prism and $\widetilde{F}$ is one of its triangular face, then the blowdown of $Q$ of the face $\widetilde{F}$ on one of its vertex is a tetrahedron. 

Now we investigate how the blowdown of a simple polytope affects its retraction sequences if the polytope remains simple after the blowdown.

\begin{theorem}\label{retraction sequence of blowdown}
Let $Q'$ be the blowdown of an $n$-dimensional simple polytope $Q$ of the facet $\widetilde{F}$ on $F$, and $Q'$ is simple.
 For a retraction sequence $\{(B_\ell, E_\ell, b_\ell)\}_{\ell=1}^m$ of $Q$ where $m=|V(Q)|$, there exists a retraction sequence $\{B_t', E_t', b_t'\}_{t=1}^{m-k}$ of $Q'$ which preserves the ordering on vertices. 
\end{theorem}

\begin{proof}
We adhere to the notations from the proof of Lemma \ref{Q' is simple}. 
Also, recall $F$ is a face of dimension $d$ in $Q$. Then from the converse part of Lemma \ref{Q' is simple}, $\widetilde{F}$ is homeomorphic to $F \times \Delta^{n-d-1}$ as a manifold with corners.
Let $V(F):=\{b_{\ell_1}, \dots, b_{\ell_k}\} \subset V(Q)$ be the vertices of $F$ such that $\ell_1 < \dots <\ell_k$.  
% First let $b_1 \not \in V(\bar{F})$. We choose $\tilde{f}(b_1)$ as $b'_1$. Otherwise if $b_1 \in V(\bar{F})$, then it must of the form $b_{\ell_1} \times \{0\}$ or $b_{\ell_1} \times \{1\}$. In this case 
 %We take $b'_1=\tilde{f}(b_1)$ where $\tilde{f}$ is defined in (\ref{define tilde f}).
%Being a vertex in an $n$-dimensional simple polytope, $b'_1$ has a local neighborhood homeomorphic to $\mathbb{R}_{\geq 0}^n$ as a manifold with corners. So $B'_1=E'_1=Q'$  gives the first entry of the sequence $(B'_1,E'_1,b'_1)$.
We construct a retraction sequence $\{(B'_t, E'_t, b'_t)\}_{t=1}^{m-k}$ of $Q'$ inductively. First we define $(B'_1, E'_1, b'_1) := (Q', Q', b'_1)$ where $b'_1 := \widetilde{f} (b_1)$, since $Q'$ is a simple polytope. Now we may encounter the following 3 cases to construct the second triple for a retraction sequence of $Q'$.

{\bf{Case 1 of the 2-nd step:}}
Let $b_2$ be neither in $V(\widetilde{F})$ nor adjacent to any vertex in $V(\widetilde{F})$ and $C_1' := \cup \{E:E~\text{is a face of}~B'_1~ \text{containing the vertex}~ b'_1\}$. Then we take $b'_2:=\tilde{f}(b_2)$ and define 
\begin{equation}\label{blowdown b2}
B'_2 :=B'_1 \setminus C_1', \quad \mbox{and} \quad
E'_2 :=\tilde{f}(E_2).
\end{equation}
The definition of blowdown implies that $Q'$ does not change locally at the points which are not in $V(\widetilde{F})$ or adjacent to a vertex in $V(\widetilde{F})$, see Remark \ref{blowdown of face}. 
 This implies $E'_2 \cong E_2$ as a manifold with corners, in which $b'_2$ has a neighborhood homeomorphic to $\mathbb{R}_{\geq 0}^{n-1}$ as a manifold with corners. So we can construct the next triple $(B'_2,E'_2,b'_2)$.

{\bf{Case 2 of the 2-nd step:}} Let $b_2 \in V(\widetilde{F})$. If $\tilde{f}(b_2)=\tilde{f}(b_1)$, then $\tilde{f}(b_2)$ is already retracted. Then to define $b'_2$ we need to go to the next vertex $\tilde{f}(b_3)$. Otherwise, we take $b'_2:=\tilde{f}(b_2)$. So we can get $(B'_2, E'_2, b'_2)$ where $B'_2$ and $E'_2$ is defined as in (\ref{blowdown b2}).
 As $b'_2$ is connected to $b'_1$ through an edge, $b'_2$ has a neighborhood homeomorphic to $\mathbb{R}_{\geq 0}^{n-1}$ in $B'_2$ as a manifold with corners. Thus we get the second entry $(B'_2,E'_2,b'_2)$ of a retraction sequence for $Q'$.
 
 {\bf{Case 3 of the 2-nd step:}} Let $b_2$ be adjacent to a vertex in $V(\widetilde{F})$. We define $b'_2:=\tilde{f}(b_2)$ along with $B'_2$ and $E'_2$ as in (\ref{blowdown b2}).
%Similar to the \textbf{Case 2 of 2-nd step}, 
As $b'_2$ is connected to $b'_1$ by an edge, $b'_2$ has a neighborhood homeomorphic to $\mathbb{R}_{\geq 0}^{n-1}$ in $B'_2$ as a manifold with corners. Thus we get the second triple $(B'_2,E'_2,b'_2)$ for the retraction sequence of $Q'$.

%If $b_2 \not \in V(\bar{F})$ but adjacent to a vertex in $V(\bar{F})$ we take $b'_2=\tilde{f}(b_2)$. Letting $B'_2$ and $E'_2$ as of previous case, see (\ref{blowdown b2}), gives us $(B'_2,E'_2,b'_2)$ by similar arguments to the previous case.

\begin{figure}
\begin{tikzpicture}[scale=.35]

\draw (0,3)--(2,4)--(1,5)--(-1,5)--(-2,4)--cycle;
\draw (-2,4)--(-2,1)--(0,0)--(2,1)--(2,4);
\draw (0,0)--(0,3);
\draw [dashed] (-2,1)--(-1,2)--(1,2)--(2,1);
\draw [dashed] (-1,2)--(-1,5);
\draw [dashed] (1,2)--(1,5);
\draw (0,0)--(2,1)--(2,4)--(0,3)--cycle;
\node at (0,-1) {$B_1$};
%\draw [<-] (1.2,.5)--(2.5,.5);
%\node [right] at (2.5,.4) {$F$};
%\node at (.5,1.3) {$\widetilde{F}$};
%\draw[<-] (-5.5,2.5)--(-2.5,2.5);
\draw[->] (2.5,2.5)--(3.5,2.5); 
\node at (2,4) {$\bullet$};
\node [right] at (2,4) {$b_1$};

\begin{scope}[xshift=170]
    \draw (1,5)--(-1,5)--(-2,4)--(0,3);
\draw (-2,4)--(-2,1)--(0,0)--(2,1);
\draw (0,0)--(0,3);
\draw [dashed] (-2,1)--(-1,2)--(1,2)--(2,1);
\draw [dashed] (-1,2)--(-1,5);
\draw [dashed] (1,2)--(1,5);
\draw (0,0)--(2,1);
\node at (0,-1) {$B_2$};
%\draw [<-] (1.2,.5)--(2.5,.5);
%\node [right] at (2.5,.4) {$F$};
%\node at (.5,1.3) {$\widetilde{F}$};
%\draw[<-] (-5.5,2.5)--(-2.5,2.5);
\draw[->] (2.5,2.5)--(3.5,2.5); 
\node at (1,5) {$\bullet$};
\node [right] at (1,5) {$b_2$};
\end{scope}

\begin{scope}[xshift=340]
\draw (-1,5)--(-2,4)--(0,3);
\draw (-2,4)--(-2,1)--(0,0)--(2,1);
\draw (0,0)--(0,3);
\draw [dashed] (-2,1)--(-1,2)--(1,2)--(2,1);
\draw [dashed] (-1,2)--(-1,5);
\draw (0,0)--(2,1);
\node at (0,-1) {$B_3$};
%\draw [<-] (1.2,.5)--(2.5,.5);
%\node [right] at (2.5,.4) {$F$};
%\node at (.5,1.3) {$\widetilde{F}$};
%\draw[<-] (-5.5,2.5)--(-2.5,2.5);
\draw[->] (2.5,2.5)--(3.5,2.5); 
\node at (2,1) {$\bullet$};
\node [right] at (2,1) {$b_3$};
\end{scope}

\begin{scope}[xshift=510]
\draw (-1,5)--(-2,4)--(0,3);
\draw (-2,4)--(-2,1)--(0,0);
\draw (0,0)--(0,3);
\draw [dashed] (-2,1)--(-1,2)--(1,2);
\draw [dashed] (-1,2)--(-1,5);
\node at (0,-1) {$B_4$};
%\draw [<-] (1.2,.5)--(2.5,.5);
%\node [right] at (2.5,.4) {$F$};
%\node at (.5,1.3) {$\widetilde{F}$};
%\draw[<-] (-5.5,2.5)--(-2.5,2.5);
\draw[->] (2,2.5)--(3,2.5); 
\node at (1,2) {$\bullet$};
\node [above] at (1,2) {$b_4$};
\end{scope}

\begin{scope}[xshift=660]
\draw (-1,5)--(-2,4)--(0,3);
\draw (-2,4)--(-2,1)--(0,0);
\draw (0,0)--(0,3);
\draw [dashed] (-2,1)--(-1,2);
\draw [dashed] (-1,2)--(-1,5);
\node at (0,-1) {$B_5$};
%\draw [<-] (1.2,.5)--(2.5,.5);
%\node [right] at (2.5,.4) {$F$};
%\node at (.5,1.3) {$\widetilde{F}$};
%\draw[<-] (-5.5,2.5)--(-2.5,2.5);
\draw[->] (.75,2.5)--(1.75,2.5); 
\node at (0,3) {$\bullet$};
\node [above] at (0,3) {$b_5$};
\end{scope}

\begin{scope}[xshift=780]
\draw (-1,5)--(-2,4);
\draw (-2,4)--(-2,1)--(0,0);
\draw [dashed] (-2,1)--(-1,2);
\draw [dashed] (-1,2)--(-1,5);
\node at (-.5,-1) {$B_6$};
%\draw [<-] (1.2,.5)--(2.5,.5);
%\node [right] at (2.5,.4) {$F$};
%\node at (.5,1.3) {$\widetilde{F}$};
%\draw[<-] (-5.5,2.5)--(-2.5,2.5);
\draw[->] (.5,2.5)--(1.5,2.5); 
\node at (0,0) {$\bullet$};
\node [above] at (0,0) {$b_6$};
\end{scope}

\begin{scope}[xshift=890]
\draw (-1,5)--(-2,4);
\draw (-2,4)--(-2,1);
\draw [dashed] (-2,1)--(-1,2);
\draw [dashed] (-1,2)--(-1,5);
\node at (-.75,-1) {$B_7$};
%\draw [<-] (1.2,.5)--(2.5,.5);
%\node [right] at (2.5,.4) {$F$};
%\node at (.5,1.3) {$\widetilde{F}$};
%\draw[<-] (-5.5,2.5)--(-2.5,2.5);
\draw[->] (0,2.5)--(1,2.5); 
\node at (-1,5) {$\bullet$};
\node [above] at (-1,5) {$b_7$};
\end{scope}

\begin{scope}[xshift=990]
\draw (-2,4)--(-2,1);
\draw [dashed] (-2,1)--(-1,2);
\node at (-.75,-1) {$B_8$};
%\draw [<-] (1.2,.5)--(2.5,.5);
%\node [right] at (2.5,.4) {$F$};
%\node at (.5,1.3) {$\widetilde{F}$};
%\draw[<-] (-5.5,2.5)--(-2.5,2.5);
\draw[->] (0,2.5)--(1,2.5); 
\node at (-1,2) {$\bullet$};
\node [above] at (-1,2) {$b_8$};
\end{scope}

\begin{scope}[xshift=1090]
\draw (-2,4)--(-2,1);
\node at (-1,-1) {$B_9$};
%\draw [<-] (1.2,.5)--(2.5,.5);
%\node [right] at (2.5,.4) {$F$};
%\node at (.5,1.3) {$\widetilde{F}$};
%\draw[<-] (-5.5,2.5)--(-2.5,2.5);
\draw[->] (-1,2.5)--(.0,2.5); 
\node at (-2,4) {$\bullet$};
\node [above] at (-2,4) {$b_9$};
\end{scope}

\begin{scope}[xshift=1190]
\node at (-2,1) {$\bullet$}; 
\node [above] at (-2,1) {$b_{10}$};
\node at (-1,-1) {$B_{10}$};
\end{scope}

\begin{scope}[yshift=-250]
\draw (1,5)--(-1,5)--(-2,4);
\draw (-2,4)--(-2,1)--(0,0)--(2,1);
\draw [dashed] (-2,1)--(-1,2)--(1,2)--(2,1);
\draw [dashed] (-1,2)--(-1,5);
\draw [dashed] (1,2)--(1,5);
\draw (-2,4)--(0,0);
\draw (1,5)--(2,1);
%\draw [very thick] (0,0)--(2,1);
%\draw [<-] (1.2,.5)--(2.5,.5);
%\node [right] at (2.5,.4) {$F$};
\node at (0,-1) {$B'_1$};
\node at (2,1) {$\bullet$};
\node [right] at (2,1) {$b'_1$};
\draw[->] (2.5,2.5)--(3.5,2.5); 
\end{scope}

\begin{scope}[yshift=-250, xshift=170]
\draw (1,5)--(-1,5)--(-2,4);
\draw (-2,4)--(-2,1)--(0,0);
\draw [dashed] (-2,1)--(-1,2)--(1,2);
\draw [dashed] (-1,2)--(-1,5);
\draw [dashed] (1,2)--(1,5);
\draw (-2,4)--(0,0);
%\draw [very thick] (0,0)--(2,1);
%\draw [<-] (1.2,.5)--(2.5,.5);
%\node [right] at (2.5,.4) {$F$};
\node at (0,-1) {$B'_2$};
\node at (1,5) {$\bullet$};
\node [right] at (1,5) {$b'_2$};
\draw[->] (2.5,2.5)--(3.5,2.5); 
\end{scope}

\begin{scope}[yshift=-250, xshift=340]
\draw (-1,5)--(-2,4);
\draw (-2,4)--(-2,1)--(0,0);
\draw [dashed] (-2,1)--(-1,2)--(1,2);
\draw [dashed] (-1,2)--(-1,5);
\draw (-2,4)--(0,0);
%\draw [very thick] (0,0)--(2,1);
%\draw [<-] (1.2,.5)--(2.5,.5);
%\node [right] at (2.5,.4) {$F$};
\node at (0,-1) {$B'_3$};
\node at (1,2) {$\bullet$};
\node [above] at (1,2) {$b'_3$};
\draw[->] (2.5,2.5)--(3.5,2.5); 
\end{scope}

\begin{scope}[yshift=-250, xshift=510]
\draw (-1,5)--(-2,4);
\draw (-2,4)--(-2,1)--(0,0);
\draw [dashed] (-2,1)--(-1,2);
\draw [dashed] (-1,2)--(-1,5);
\draw (-2,4)--(0,0);
%\draw [very thick] (0,0)--(2,1);
%\draw [<-] (1.2,.5)--(2.5,.5);
%\node [right] at (2.5,.4) {$F$};
\node at (0,-1) {$B'_4$};
\node at (0,0) {$\bullet$};
\node [above] at (0.3,0) {$b'_4$};
\draw[->] (2,2.5)--(3,2.5); 
\end{scope}

\begin{scope}[yshift=-250, xshift=680]
\draw (-1,5)--(-2,4);
\draw (-2,4)--(-2,1);
\draw [dashed] (-2,1)--(-1,2);
\draw [dashed] (-1,2)--(-1,5);
%\draw [very thick] (0,0)--(2,1);
%\draw [<-] (1.2,.5)--(2.5,.5);
%\node [right] at (2.5,.4) {$F$};
\node at (-.5,-1) {$B'_5$};
\node at (-1,5) {$\bullet$};
\node [right] at (-1,5) {$b'_5$};
\draw[->] (1.5,2.5)--(2.5,2.5); 
\end{scope}

\begin{scope}[yshift=-250, xshift=840]
\draw (-2,4)--(-2,1);
\draw [dashed] (-2,1)--(-1,2);
%\draw [very thick] (0,0)--(2,1);
%\draw [<-] (1.2,.5)--(2.5,.5);
%\node [right] at (2.5,.4) {$F$};
\node at (-.75,-1) {$B'_6$};
\node at (-1,2) {$\bullet$};
\node [right] at (-1,2) {$b'_6$};
\draw[->] (1.5,2.5)--(2.5,2.5); 
\end{scope}

\begin{scope}[yshift=-250, xshift=980]
\draw (-2,4)--(-2,1);
%\draw [very thick] (0,0)--(2,1);
%\draw [<-] (1.2,.5)--(2.5,.5);
%\node [right] at (2.5,.4) {$F$};
\node at (-.75,-1) {$B'_7$};
\node at (-2,4) {$\bullet$};
\node [right] at (-2,4) {$b'_7$};
\draw[->] (1.2,2.5)--(2.2,2.5); 
\end{scope}

\begin{scope}[yshift=-250, xshift=1140]
%\draw [very thick] (0,0)--(2,1);
%\draw [<-] (1.2,.5)--(2.5,.5);
%\node [right] at (2.5,.4) {$F$};
\node at (-1,-1) {$B'_8$};
\node at (-2,1) {$\bullet$};
\node [above] at (-2,1) {$b'_8$};
\end{scope}

\end{tikzpicture}
\caption{Induced retraction sequence on a blowdown.}
\end{figure}

Continuing a similar way, suppose that we are at $t$-th step of a retraction of $Q'$. In the meantime, we are at $i_t$-th step of retraction of $Q$ where $t \leq i_t$. At this step, three cases may arise.

{\bf{Case 1 of the $t$-th step:}} Let $b_{i_t}$ is neither in $V(\widetilde{F})$ nor adjacent to any vertex in $V(\widetilde{F})$ and $C_{t-1}':= \cup \{E:E~ \text{is a face of} ~B'_{t-1}~\text{containing the vertex}~ b'_{t-1}\}$. Then define
\begin{equation}
b'_t:=\tilde{f}(b_{i_t}), ~~
 B'_t:=B'_{t-1} \setminus C_{t-1}', ~~ \mbox{and} ~~ 
 E'_t:=\tilde{f}(E_{i_t}) (\cong E_{i_t} ~\text{as a manifold with corners}).
\end{equation}
This gives $t$-th triple $(B'_t, E'_t, b'_t)$ of a retraction sequence for $Q'$. 

{\bf{Case 2 of the $t$-th step:}} Let $b_{i_t} \in V(\widetilde{F})$.
% then it is of the form $b_{\ell_j} \times \{0\}$ or $b_{\ell_j} \times \{1\}$ for some $j \in \{1, \dots ,k\}$. 
If $\widetilde{f}(b_{i_t}) \in \{b'_1, \dots, b'_{t-1}\}$ then $\widetilde{f}(b_{i_t})$ is already retracted. In this situation, we need to move to the next vertex in the sequence $\{b_1,\dots,b_m\}$ to get $b'_{t}$ in $Q'$.
Otherwise, we define 
\begin{equation}\label{blowdown bt}
b'_t:=\tilde{f}(b_{i_t}), ~~
 B'_t:=B'_{t-1} \setminus C_{t-1}', ~~ \mbox{and} ~~ 
E'_t:=\tilde{f}(E_{i_t}) \text{ or } E'_t \text{ is a face of } \tilde{f}(E_{i_t}) \text{ containing } b'_t,
\end{equation}
where $C_{t-1}':= \cup \{E:E~ \text{is a face of} ~B'_{t-1}~\text{containing the vertex}~ b'_{t-1}\}$.
If $b_{i_t}$ is connected to some vertex $b_u$ in $Q$ %where $u >i_t$ 
and $\tilde{f}(b_u)=b'_v$, then $$\dim(E'_t)= \dim(E_{i_t})- \#\{b_u \in V(Q) ~| ~ %\tilde{f}(b_u)=b'_v \text{ with }
u>i_t \text{ and } v<t\}.$$
%with $v<t$, then $\dim(E'_t)$ is less than $\dim(E_{i_t})$ by the number of such $b_u \in V(Q)$.
Note that from the ordering of the vertices in the retraction sequence of $Q$, we have at least one $b_u \in V(Q)$ connected to $b_{i_t}$ such that $u>i_t$ except when $b_{i_t}=b_m$ is the last vertex in the retraction sequence of $Q$. 
Also note that since the vertices $b_u \in V(Q)$ connected to $b_{i_t}$ with $u<i_t$ are retracted before $b_{i_t}$, they does not affect $\dim(E_{i_t})$ and eventually $\dim(E'_t)$.
Thus \eqref{blowdown bt} provides the $t$-th entry of a retraction sequence of $Q'$ for this case.

{\bf{Case 3 of the $t$-th step:}} Let $b_{i_t}$ is adjacent to a vertex in $V(\widetilde{F})$.
% then it is of the form $b_{\ell_j} \times \{0\}$ or $b_{\ell_j} \times \{1\}$ for some $j \in \{1, \dots ,k\}$. 
 We define $$b'_t:=\tilde{f}(b_{i_t}), ~B'_t~ \mbox{as in} \eqref{blowdown bt}
~\mbox{and} ~E'_t:=\tilde{f}(E_{i_t}) ~\mbox{or} ~E'_t ~\mbox{is a face of} \tilde{f}(E_{i_t}) ~\mbox{containing} ~b'_t.$$
An argument similar to the previous case gives $t$-th triplet $(B'_t,E'_t,b'_t)$ for this case.

%If $b_i$ is connected to some vertex $b_u$ in $Q$ where $u >i$ but $\tilde{f}(b_u)=b'_v$ with $v<t$, then dim$(E'_t)$ is less than dim$(E_i)$ by the number of such $b_u \in Q$.
%This provides $t$-th entry of a retraction sequence of $Q'$.
% If $b_i \not \in V(\bar{F})$ but adjacent to a vertex in $V(\bar{F})$, then by taking $f(b_i)=b'_t$ and defining $B'_t$ and $E'_t$ as in (\ref{blowdown bt}) we get $(B'_t, E'_t, b'_t)$ by similar arguments.

 Therefore, by the inductive process, we get $\{(B'_t, E'_t, b'_t)\}_{t=1}^{m-k}$ as the desired induced retraction sequence of $Q'$ from $Q$.
 
% We note that $Q'$ may possess different retraction sequences.
\end{proof}
 The retraction sequence $\{(B'_t, E'_t, b'_t)\}_{t=1}^{m-k}$ is called an induced retraction sequence of $Q'$.

We remark that though the blowdown $Q'$ in Figure \ref{Fig_blowdown but not simple} is not a simple polytope, it may induce a retraction sequence from that of $Q$ if the retraction sequence starts with the vertices in the base.
 On the other hand, if the retraction sequence of $Q$ in Figure \ref{Fig_blowdown but not simple} starts with any vertex of the top square, then $Q'$ doesn't induce a retraction sequence from $Q$. Thus we need $Q'$ to be simple in Theorem \ref{retraction sequence of blowdown}.

\begin{corollary}\label{cor_ret_k_wed}
Let $\{(B_\ell, E_\ell, b_\ell)\}_{\ell=1}^m$ be a  retraction sequence of $Q$, where $m=|V(Q)|$. There always exists a retraction sequence $\{(B_t, E_t, b_t)\}_{t=1}^u$ for a polytopal $k$-wedge $Q_F(k)$ of $Q$ at $F$, where $u=(k+1)m-k\alpha$ and $\alpha=|V(F)|$.
\end{corollary}
\begin{proof}
Proposition \ref{retraction of F times Delta} gives an induced retraction sequence $\{(\widebar{B}_j, \widebar{E}_j, \widebar{b}_j)\}_{j=1}^{(k+1)m}$ on $Q \times \Delta^k$ such that $\widebar{E}_{(k+1)\ell - (k+1-s)} = E_{\ell} \times \Delta^{k+1-s}$ for $\ell =1, ... , m$ and $s=1, ..., k+1$. Then the claim follows from Corollary \ref{cor_gen of poly wedge cons} and Theorem \ref{retraction sequence of blowdown}. 
\end{proof}

 Moreover, we get the following if $\{(\Delta^{k+1-s}, \Delta^{k+1-s}, e_s)\}_{s=1}^{k+1}$ is a retraction sequence of $\Delta^k$.

%the vertex of $\Delta^{\beta}$ not contained in the facet $\Delta^{\beta - 1}$ of $\Delta^{\beta}$ for $\beta = 0, \dots, k$. 
\begin{corollary}\label{cor_k wed_ret}
Let $Q$ be a simple polytope with a facet $F$ such that $|V(F)|=\alpha$ and there exists a retraction sequence $\{(B_{\ell}, E_{\ell}, v_{\ell})\}_{\ell=1}^m$ such that the vertices of $F$ to be retracted at the end. For $u=(k+1)(m-\alpha)+\alpha$ there exists a retraction sequence $\{(B'_t, E'_t, b'_t)\}_{t=1}^u$ for $Q_F(k)$ such that 
\begin{enumerate}
    \item $E'_{(k+1)\ell - (k+1-s)} = E_{\ell} \times \Delta^{k+1-s}$ for $\ell =1, \dots , m-\alpha$ and $s=1, \dots, k+1$, and $E'_t=E_{m- \alpha +\ell}$
for $t=(k+1)(m-\alpha)+\ell$ and $\ell=1, \dots, \alpha$. 

\item $b'_{(k+1)\ell - (k+1-s)} = (v_{\ell}, e_{s})$ for $\ell =1, \dots , m-\alpha$ and $s=1, \dots, k+1$, and $b'_t=v_{m- \alpha +\ell}$
for $t=(k+1)(m-\alpha)+\ell$ and $\ell=1, \dots, \alpha$.
\end{enumerate}

\end{corollary}

%\begin{remark}\label{blowdown of a simple polytope may not be a simple polytope}
%Blowdown of a simple polytope may not be a simple polytope as you can see in Figure \ref{Fig_blowdown but not simple}.
% Please note that in Figure \ref{Fig_blowdown but not simple} if we consider the retraction sequence of $Q$ in some particular pattern then it induces a retraction sequence of $Q'$.  
% \end{remark}

\section{Blowdowns of quasitoric orbifolds and torsions in their integral cohomologies}
\label{Sec_Blowdown of quasitoric orbifolds}

In this section, we study the effects of blowdowns of simple polytopes on their corresponding quasitoric orbifolds.
 Note that the blowdown of a quasitoric orbifold is discussed in \cite{GP2}. However, we study blowdowns of quasitoric orbifolds in more generality.
We also investigate the torsions in the integral cohomology of quasitoric orbifolds after blowdowns and prove no new torsion arises in certain blowdowns.
 We adhere to the notation of previous sections.

Let $Q$ be an $n$-dimensional simple polytope with $\mathcal{F}(Q)=\{\widetilde{F},F_1,\dots,F_r\}$.
 Consider a blowdown $Q'$ of $Q$ of the facet $\widetilde{F}$ on a face $F$ such that $Q'$ is a simple polytope. Let $\mathcal{F}(Q')=\{F'_1,\dots,F'_r\}$ as in Remark \ref{define f}. 
Let $\lambda \colon \mathcal{F}(Q) \to \ZZ^n$ and $\lambda' \colon \mathcal{F}(Q') \to \ZZ^n$ be two  $\mathcal{R}$-characteristic functions on $Q$ and $Q'$ respectively such that
\begin{equation}\label{defn of lambda'}
\lambda'(F'_i)=\lambda'(\tilde{f}(F_i))=\lambda(F_i)
\end{equation}
for $1 \leq i \leq r$ where $\widetilde{f}$ is defined in \eqref{define tilde f}.
Then we call the pair $(Q', \lambda')$ a \emph{restriction} of $(Q,\lambda)$. %So, we have two different $2n$-dimensional quasitoric orbifolds $X(Q, \lambda)$ and $X(Q',\lambda').$ 

\begin{definition}[Blowdown of a quasitoric orbifold]\label{blowdown definition}
Let $(Q,\lambda)$ and $(Q', \lambda')$ be$\mathcal{R}$-characteristic pairs such that $(Q',\lambda')$ is a restriction of $(Q,\lambda)$. Then the quasitoric orbifold $X(Q', \lambda')$ is called a blowdown of $X(Q,\lambda)$.
% corresponding to the faces $F \subset \widetilde{F}$.
\end{definition}

\begin{example}
Let $(Q, \lambda)$ be an $\mathcal{R}$-characteristic function and $Q'$ a blowdown of $Q$ such that $Q'$ is simple. Then the natural restriction of $\lambda$ on $\mathcal{F}(Q')$ using Remark \ref{define f} and \eqref{defn of lambda'} may not be an   $\mathcal{R}$-characteristic function. 
For example, consider a blowdown of a cube as in Figure \ref{Fig_blowdown of a cube} and the $\mathcal{R}$-characteristic function $\lambda$ on the facets of $Q$ by
\begin{align}
\lambda(F_0)=(1,0,0), \quad \lambda(F_1)=(1,0,0), \quad
\lambda(F_2)=(2,3,5), \\ \lambda(F_3)=(1,3,2),\quad
\lambda(F_4)=(4,1,0), \quad \lambda(\widetilde{F})=(1,0,1).\nonumber
\end{align}
If we define $\lambda' \colon \mathcal{F}(Q') \to \ZZ^3$ by natural restriction following \eqref{defn of lambda'}, then $\lambda'$ is not an $\mathcal{R}$-characteristic function on $Q'$ since $$\det[\lambda'(F'_0),\lambda'(F'_1),\lambda'(F'_4)]=0$$ where $F'_0\cap F'_1 \cap F'_4$ is a vertex in $Q'$. Thus the pair $(Q',\lambda')$ does not determine any quasitoric orbifold.
\end{example}

This justifies our definition of restriction of an $\mathcal{R}$-characteristic pair as well as the blowdown of quasitoric orbifolds. Next, we give a sufficient condition when the natural restriction is an $\mathcal{R}$-characteristic function.

%Let $F$ be a face and $\widetilde{F}$ a facet of an $n$-dimensional simple polytope $Q$ with $F \subset \widetilde{F}$. 
Let $Q'$ be a blowdown of $Q$ of the facet $\widetilde{F}$ on $F$.
Then by Lemma \ref{Q' is simple}, $Q'$ is simple if and only if $\widetilde{F} \cong F \times \Delta^{n-d-1}$ as a manifold with corners.
If $b \in Q$ be any vertex such that $b \not \in \widetilde{F}$, then the facets adjacent to $\widetilde{f}(b)$ remain the same, see \eqref{vertex not in facet}.
Now let us fix a vertex $b \in \widetilde{F}$. So $b=(\bigcap_{j=1}^{n-1} F_{i_j}) \cap \widetilde{F}$ as in \eqref{b_l_i}.
Similar construction as in the proof of Lemma \ref{Q' is simple} leads us to $$\widetilde{f}(b)= (\bigcap_{j=1}^{n-1} \widetilde{f}(F_{i_j})) \cap \widetilde{f}(P_b),$$ for a unique facet $P_b$ of $Q$ such that $b \not \in P_b$, see \eqref{tilde f b_l_i}.
%Note that image of $V(\widetilde{F})$ through $\widetilde{f}$ is same as the image of $V(F)$. So considering $b \in V(F)$ is enough.
We define a set
$$S_{b} :=\{\lambda(F_{i_1}),\dots,\lambda(F_{i_{n-1}}),\lambda(P_b)\},$$
 for each vertex $b \in \widetilde{F}$.
As a vertex of $ \widetilde{F}$, $b$ can be considered as   $(b_F,v_q)$ for some $b_F \in V(F)$ and $v_q \in V(\Delta^{n-d-1})$. Notice that for any $v, ~v' \in V(\Delta^{n-d-1})$ we have $$S_{({b}_F,v)}=S_{(b_F,v')}.$$
So we denote $S_{b_F}:=S_{(b_F,v)}$. Then, we can conclude the following.
 
\begin{proposition}\label{S_b is linearly independent} 
Let $Q'$ be the blowdown of $Q$ of the facet $\widetilde{F}$ on $F$. If $S_{b_F}$ is linearly independent for each $b_F \in V({F})$, then the pair $(Q', \lambda')$  is an $\mathcal{R}$-characteristic pair as well as a restriction of $(Q, \lambda)$, where $\lambda'$ is defined in \eqref{defn of lambda'}. 
\end{proposition}

\begin{example}\label{blowdown of cube}
Let $Q$ be a cube as in Figure \ref{Fig_blowdown of a cube}. We define $\lambda$ on $\mathcal{F}(Q)$ by

\begin{align}\label{char vec on cube for blowdown case}
\lambda(F_0)=(0,2,1), \quad \lambda(F_1)=(1,1,2), \quad
\lambda(F_2)=(0,1,1), \\ \lambda(F_3)=(1,0,1),\quad
\lambda(F_4)=(1,0,0), \quad \lambda(\widetilde{F})=(1,3,3).\nonumber
\end{align}
This gives an $\mathcal{R}$-{characteristic pair} $(Q, \lambda)$ and consequently a quasitoric orbifold $X(Q,\lambda)$. In Figure \ref{Fig_blowdown of a cube}, $Q'$ is the blowdown of $Q$ of the face $\widetilde{F}$ on the face $F$.
 We define $\lambda'$ on the facets of $Q'$ by (\ref{defn of lambda'}). This gives a restriction $(Q',\lambda')$ of $(Q, \lambda)$.
So $X(Q',\lambda')$ is a blowdown of $X(Q, \lambda)$.\qed
% corresponding to the faces $F \subset \widetilde{F}$.

\begin{figure}
\begin{tikzpicture}[scale=.6]
\draw (0,0)--(0,2)--(-1,3)--(-1,1)--cycle;
\draw[fill=yellow, opacity=.4] (0,0)--(2,0)--(2,2)--(0,2)--cycle;
\draw [thick] (0,0)--(2,0);
\draw (2,2)--(1,3)--(-1,3);
\draw [dashed] (-1,1)--(1,1)--(1,3);
\draw [dashed] (1,1)--(2,0);

\draw (0,2.5)--(0,3.4);
\node [above] at (0,3.4) {$F_0$};

\draw (.5,0) to [out=240, in=330] (-.5,-.3);
\draw[dotted,thick] (.5,0)--(.7,.3);
\node [left] at (-.5,-.3) {$F_1$};

\draw [dotted, thick] (1.5,1.7)--(2,1.9);
\draw (2,1.9)--(2.5,2.1);
\node [right] at (2.5,2.1) {$F_2$};

\draw[dotted,thick] (-.7,1.8)--(-1,2);
\draw (-1,2)--(-1.7,2.5);
\node [left] at (-1.7,2.5) {$F_3$};

\draw (-.7,1.3)--(-1.7,1.3);
\node [left] at (-1.7,1.3) {$F_4$};

\draw (1.7,.7)--(2.3,.7);
\node [right] at (2.3,.7) {$\widetilde{F}$};

\draw [thick] (0,0)--(2,0);
\draw [->] (1.7,-.4)--(1.5,0);
\node [right] at (1.5,-.6) {$F$};

\node at (1,-1.3) {$Q$};

\draw[->] (3,1.5)--(5.5,1.5);
%\node [above] at (4,1.5) {$Blowdown$};
%\node [right] at (3,1.2) {$of~\widebar{F}~on~F$};

\begin{scope}[xshift=200]
\draw (0,0)--(-1,3)--(-1,1)--cycle;
\draw (0,0)--(2,0)--(1,3)--(-1,3);
\draw [dashed] (-1,1)--(1,1)--(1,3);
\draw [dashed] (1,1)--(2,0);

\node at (1,-1.3) {$Q'$};

\draw (0,2.5)--(0,3.4);
\node [above] at (0,3.4) {$F'_0$};

\draw (.5,0) to [out=240, in=330] (-.5,-.3);
\draw[dotted,thick] (.5,0)--(.7,.3);
\node [left] at (-.5,-.3) {$F'_1$};

\draw [dotted, thick] (1.2,1.3)--(1.5,1.5);
\draw (1.5,1.5)--(2,1.8);
\node [right] at (2,1.8) {$F'_2$};

\draw[dotted,thick] (-.7,1.8)--(-1,2);
\draw (-1,2)--(-1.7,2.5);
\node [left] at (-1.7,2.5) {$F'_3$};

\draw (-.6,1.3)--(-1.2,.5);
\node [left] at (-1.2,.5) {$F'_4$};

\draw [thick] (0,0)--(2,0);
\draw [->] (1.7,-.4)--(1.5,0);
\node [right] at (1.5,-.6) {$F'$};

\end{scope}

\end{tikzpicture}
\caption{Blowdown of a cube of $\widetilde{F}$ on $F$.}
\label{Fig_blowdown of a cube}
\end{figure} 

\end{example}

%In the remaining part of this section we are going to consider $Q$ as in Example \ref{example of blowdown of quasitoric orbifold} unless stated otherwise. 
%We consider $Q$ with a facet $\widetilde{F}$ homeomorphic to $F \times I$ as a manifold with corners where $F \subset \widetilde{F}$ and $I=[0,1]$ and $Q'$ as the blowdown of $Q$  of the face $\widetilde{F}$ on $F$.

%\begin{example}\label{Rmk_linear combination in blowdown}
%We adhere the notations from Example \ref{example of blowdown of quasitoric orbifold}.
%If we assume that $\lambda(\widetilde{F})=c_{0}\lambda(F_0)+c_1\lambda(F_1)$ for nonzero $c_0$ and $c_1$ in $\mathbb{Q}$, then the condition (\ref{blowdown char vec}) follows automatically. From the construction, 
%\small{\begin{equation}
% {\text{det}[\lambda(F_2)^t, \dots, \lambda (F_{n-1})^t, \lambda(\widetilde{F})^t,\lambda(F_0)^t]=c_1 \text{det}[\lambda(F_2)^t, \dots, \lambda (F_{n-1})^t,\lambda(F_1)^t,\lambda(F_0)^t]}
%\end{equation}}
%This implies $ \text{det}[\lambda(F_2)^t, \dots ,\lambda (F_{n-1})^t,\lambda(F_1)^t,\lambda(F_0)^t]$ is nonzero as the left side of the equation and $c_1$ are nonzero. Thus $S_v$ is linearly independent.
%\end{example}

\begin{remark}\label{Rem_Immidiate use}
Let $Q'$ be a blowdown of an $n$-dimensional simple polytope $Q$ of the facet $\widetilde{F}=F \times \Delta^{n-d-1}$ on a $d$-dimensional face $F$.
If $\mathcal{F}({\Delta^{n-d-1}})=\{F_1^\Delta, \dots, F_{n-d}^\Delta\}$, then $F \times F_s^{\Delta}$ are some facets of $\widetilde{F}$ and $F \times F_s^{\Delta}=\widetilde{F} \cap P_s$ for a unique $P_s \in \mathcal{F}(Q)$ for $1 \leq s \leq (n-d).$
\end{remark}

\begin{proposition}\label{Propn Char vec over Q'}
Let $Q'$ be a blowdown of an $n$-dimensional simple polytope $Q$ of the facet $\widetilde{F}=F \times \Delta^{n-d-1}$ on a $d$-dimensional face $F$.
%If $\mathcal{F}({\Delta^{n-d-1}})=\{F_1^\Delta, \dots, F_{n-d}^\Delta\}$, then $F \times F_s^{\Delta}$ are some facets of $\widetilde{F}$ and $F \times F_s^{\Delta}=\widetilde{F} \cap P_s$ where $P_s \in \mathcal{F}(Q)$ for $1 \leq s \leq (n-d).$
Let $\lambda$ be an $\mathcal{R}$-characteristic function on $Q$ such that 
 \begin{equation}\label{lambda after blowdown}
 \lambda(\widetilde{F})=\sum_{s=1}^{n-d} c_s \lambda(P_s),
 \end{equation}
for some $c_s \in \QQ \setminus \{0\}$, where $P_s$ are described in Remark \ref{Rem_Immidiate use} for $1 \leq s \leq (n-d)$. Then $(Q', \lambda')$ is a restriction of $(Q, \lambda)$ where $\lambda'$ is defined as \eqref{defn of lambda'}.
\end{proposition}

\begin{proof}
%Let $Q'$ be a blowdown of an $n$-dimensional simple polytope $Q$ of the facet $\widetilde{F}$ on a $d$-dimensional face $F$ such that $Q'$ is simple. Thus by converse part of Lemma \ref{Q' is simple}, $\widetilde{F}$ is homeomorphic to $F \times \Delta^{n-d-1}$ as a manifold with corners.

%Let $\{F_1^\Delta, \dots, F_{n-d}^\Delta\}$ be the facets of $\Delta^{n-d-1}$. Then $F \times F_s^\Delta$ for $1 \leq s \leq (n-d)$ are some facets of $\widetilde{F}$ as well as they are codimension $2$ faces in $Q$. Thus $F \times F_s^\Delta = \widetilde{F} \cap P_s$ for some unique facets $P_s$ in $Q$ for $1 \leq s \leq (n-d)$.

Let $b \in V(F)$. So $\widetilde{f}(b) \in V(F') \subset V(Q')$.
The arguments in the proof of Lemma \ref{Q' is simple} and \eqref{b_l_i}, \eqref{tilde f b_l_i} give us
\begin{align}
b = (\bigcap_{j=1}^{n-1} F_{i_j}) \cap \widetilde{F} \quad \text{ and } \quad
\widetilde{f}(b) = (\bigcap_{j=1}^{n-1} \widetilde{f}(F_{i_j})) \cap \widetilde{f}(P_b),
\end{align}
for some unique facets $F_{i_1},\dots,F_{i_{n-1}}, \widetilde{F}$ and $P_b$ of $Q$. Note that $P_b$ is the unique facet in $\{P_s \colon 1 \leq s \leq (n-d)\}$ such that it does not contain the vertex $b$.
If we define $\lambda' \colon \mathcal{F}(Q') \to \ZZ^n$ by using \eqref{defn of lambda'} from $\lambda$ on $Q$ then
\begin{equation}\label{det at b_l_i}
\det[\lambda(F_{i_1}), \dots, \lambda(F_{i_{n-1}}), \lambda(\widetilde{F})]=c_b \det[\lambda'(\widetilde{f}(F_{i_1})), \dots,\lambda'(\widetilde{f}(F_{i_{n-1}})), \lambda'(\widetilde{f}(P_b))],
\end{equation}
where $c_b \in \{c_1, \dots, c_{n-d}\}$ is the coefficient of $\lambda(P_b)$ in \eqref{lambda after blowdown}.
This implies the vectors in $S_b$ are linearly independent, that is the vectors assigned to facets adjacent to the vertex $\widetilde{f}(b)$ in $Q'$ are linearly independent. 
We can do the above construction for each vertex in $V(F)$.
Thus at each vertex of $F'$ in $Q'$ the vectors assigned to the adjacent facets are linearly independent. 

Now let $b' \in V(Q') \setminus V(F')$. Then there exists unique $b \in V(Q) \setminus V(\widetilde{F})$ such that $\widetilde{f}(b)=b'$. The vectors assigned to the adjacent facets of $b $ in $Q$ and $b'$ in $Q'$ are the same. 
So the induced $\lambda'$ using \eqref{defn of lambda'} becomes an $\mathcal{R}$-characteristic function on $Q'$. Thus $(Q',\lambda')$ is a restriction of $(Q,\lambda)$. % and $X(Q',\lambda')$ is a blowdown of $X(Q,\lambda)$. 
\end{proof}

%{\color{blue} 
Let $E$ be a face of $Q$ such that $E\cap \widetilde{F}=E^F\times \Delta^{q}$ for some $q>0$ and $E\cap\widetilde{F}$ is a facet of $E$. Then using Lemma \ref{lem:blowdown of face}, $\widetilde{f}(E)$ is a blowdown of $E$ of the facet $E\cap \widetilde{F}$ to $E^F$. The next lemma deduces that if the $\mathcal{R}$-characteristic function $\lambda$ on $Q$ satisfies \eqref{lambda after blowdown}, then a similar relation also holds for $\lambda_E$. Recall the facet set $\mathcal{F}(\Delta^{n-d-1})=\{F_1^\Delta, F_2^\Delta,\dots ,F_{n-d}^\Delta\}$ from Remark \ref{Rem_Immidiate use}. Let $\{F_{i_1}^\Delta,F_{i_2}^\Delta,\dots,F_{i_{q+1}}^\Delta\}\subset \{F_1^\Delta,F_2^\Delta,\dots,F_{n-d}^\Delta\}$ such that $F_{i_j}^\Delta\cap\Delta^q$ is a facet of $\Delta^q$ for all $j=1,2,\dots,q+1$ and $\bigcap_{j=q+2}^{n-d}F_{i_j}^{\Delta}=\Delta^q$. Thus $\{P_{i_1},P_{i_2},\dots,P_{i_{q+1}}\}\subset \{P_1,P_2,\dots,P_{n-d}\}$ such that %$\bigcap_{j=q+2}^{n-d}P_{i_j}=\Delta^q$ and 
$\{ P_{i_1} \cap E,\dots,P_{i_{q+1}} \cap E \} \subset \mathcal{F}(E)$.

\begin{lemma}\label{projection lemma}
    If $\lambda$ be an $\mathcal{R}$-characteristic function on $Q$ satisfying \eqref{lambda after blowdown} then $$\lambda_{E}(E\cap \widetilde{F})=(\sum_{j=1}^{q+1} c_{i_j} d_j \lambda_E (E\cap P_{i_j}))/d_E,$$ for some positive integers $d_E$ and $d_j$'s where $j=1,2, \dots,q+1$.
\end{lemma}

\begin{proof}
     The projection map is defined by
$$\rho_E:\mathbb{Z}^n \to \mathbb{Z}^{\dim(E)}.$$
This map is $\mathbb{Z}$-linear and any $\mathbb{Z}$-linear map is $\mathbb{Q}$-linear.
From the definition of the induced $\mathcal{R}$-{characteristic function} on $E$
\begin{equation}
\lambda_E: \mathcal{F}(E) \to \mathbb{Z}^{\dim(E)},
\end{equation}
 we have 
 \begin{equation}\label{Eq_define d_i}
 \lambda_E(E \cap {P_{i_j}})=\text{prim}\{ \rho_E (\lambda(P_{i_j})) \}=\frac{\rho_E (\lambda(P_{i_j}))} {d_j}
 \end{equation}
 Thus 
\begin{align}\label{define d_F}
\lambda_E(E \cap \widetilde{F}) &= \text{prim}\{\rho_E (\lambda(\widetilde{F}))\}=\text{prim}\{\rho_E (\sum_{i=1}^{n-d} c_i\lambda(P_{i}))\}=\text{prim}\{\sum_{i=1}^{n-d} c_i\rho_E (\lambda(P_{i}))\} \\
  & =\text{prim}\{\sum_{j=1}^{q+1} c_{i_j} d_j \lambda_E(P_{i_j} \cap E)\} =(\sum_{j=1}^{q+1} c_{i_j} d_j \lambda_E(P_{i_j} \cap E))/d_E \nonumber,
\end{align}
for some unique positive integer $d_E$.
\end{proof}

Let $(Q,\lambda)$ be an $\mathcal{R}$-characteristic pair and $\{(B_\ell,E_\ell,b_\ell)\}_{\ell=1}^m$ be a retraction sequence of $Q$.
Then we denote $|G_{B_\ell}(b_\ell)|:=|G_{E_\ell}(b_\ell)|$ for all $\ell=1, \dots, m$.

\begin{proposition}[Proposition 4.5,\cite{BSSau21}]\label{order of group}
Let $(Q, \lambda)$ be an $\mathcal{R}$-characteristic pair and $(F,\lambda_F)$ the induced $\mathcal{R}$-characteristic pair on a face $F$ with $V(F)=\{b_{\ell_1}, \dots, b_{\ell_k}\}$. If $\{(B_i^F, E_i^F, b_i^F)\}_{i=1}^k$ is an induced retraction sequence of $F$ from $\{(B_\ell, E_\ell, b_\ell)\}_{\ell=1}^m$ of $Q$ such that $B_{i}^F=B_{\ell_i} \cap F$, $b_i^F=b_{\ell_i}$ and $E_i^F$ is the maximal dimensional face of $B_i^F$ containing $b_i^F$; 
then $|G_{E_i^F}(b_i^F)|$ divides $|G_{E_{\ell_i}}(b_{\ell_i})|$. 
\end{proposition}

Next, we discuss how the singularities are affected after certain blowdowns of quasitoric orbifolds.
\begin{theorem}\label{blowdown theorem}
Let $(Q, \lambda)$ be an $\mathcal{R}$-characteristic pair satisfying the hypothesis in Proposition \ref{Propn Char vec over Q'} and  
$Q'$ a blowdown of $Q$ of the facet $\widetilde{F}$ on $F$ with  $|V(F)|=k, \dim F =d$. Let $p$ be a prime such that the following holds:
\begin{itemize}
\item[($A_1$)] There exists a retraction sequence $\{(B_i, E_i, b_i)\}_{i=1}^m$ of $Q$ such that $\gcd(|G_{E_i}(b_i)|,p)=1$ for all $i=1, \dots, m$.
\item[($A_2$)] The map $\lambda \colon \mathcal{F}(Q) \to \ZZ^n$ satisfies \eqref{lambda after blowdown} such that $\gcd(\text{\emph{denominator~of}}~c_s,p)=1$ for $s=1,\dots, (n-d)$
\item[($A_3$)] If $E_i\cap \widetilde{F}$ is a facet of $E_i$ such that $\widetilde{f}(E_i)$ is a blowdown of $E_i$ of the facet $E_i\cap \widetilde{F}$ to some face $E_i^F$ then $\gcd(d_{E_i}, p)=1$, where $\lambda_{E_i}(E_i\cap \widetilde{F})=d_{E_i}\cdot{\rho_{E_i}(\lambda(\widetilde{F}))}$. 
\end{itemize}
Then $X(Q',\lambda')$ is a blowdown of $X(Q,\lambda)$ where $\lambda'$ is defined in \eqref{defn of lambda'} and $Q'$ 
has a retraction sequence $\{{B}'_t,{E}'_t,{b}'_t\}_{t=1}^{m-k}$ such that $\gcd(|G_{{E}'_t}({b}'_t)|,p)=1$ for all $t=1, \dots,m-k$.
\end{theorem}
\begin{proof}
Let $X(Q,\lambda)$ be a quasitoric orbifold over an $n$-dimensional simple polytope $Q$ having a retraction sequence $\{(B_\ell, E_\ell, b_\ell)\}_{\ell=1}^m$.

Let $Q'$ be a blowdown of $Q$ of the facet $\widetilde{F}$ on the face $F$ such that $Q'$ is simple. Then there is a retraction sequence $\{B'_j,{E'}_j,{b'}_j\}_{j=1}^{m-k}$ of $Q'$ where $|V(F)|=k$, see the proof of Theorem \ref{retraction sequence of blowdown}.

Suppose that $\lambda$ satisfies \eqref{lambda after blowdown} and $p$ is a prime number such that $$\gcd(\text{denominator~of}~c_s,p)=1 \text{ for } s=1,\dots, (n-d).$$
 Then $X(Q',\lambda')$ is a blowdown of $X(Q, \lambda)$ by Proposition \ref{Propn Char vec over Q'}.
 
For an arbitrary vertex $b \in V(F)$, there exists $b'=\widetilde{f}(b) \in V(F')$ during the induced retraction as in the proof of Theorem \ref{retraction sequence of blowdown}. Then from \eqref{det at b_l_i}, we have
\begin{equation}\label{equation for blowdown}
|G_Q(b)|=c_s |G_{Q'}(b')|
\end{equation}
for some $s \in \{1 ,\dots,(n-d)\}$.
For the quasitoric orbifold $X(Q,\lambda)$, let us assume $$\gcd(|G_{E_i}(b_i)|,p)=1 \text{ for } i=1, \dots,m.$$ Now we want to see how the orders of the singularities $G_{E'_t}(b'_t)$ behave due to blowdown where $b'_t=\widetilde{f}(b_{i_t})$ for some $b_{i_t} \in V(Q)$.
Depending on $b_{i_t} \in V(Q)$ three cases may arise during the induced retraction of $Q'$.

{\bf{Case 1:}} Let the vertex $b_{i_t} \in V(Q)$ be neither in $V(\widetilde{F})$ nor adjacent to any vertex in $V(\widetilde{F})$. Then $E'_t$ is homeomorphic to $E_{i_t}$ as a manifold with corners, see the proof of Theorem \ref{retraction sequence of blowdown}. Thus $$|G_{E'_t}(b'_t)|=|G_{E_{i_t}}(b_{i_t})|.$$ This implies $\gcd(|G_{E'_t}(b'_t)|,p)=1$ for the vertices considered in this case.

{\bf{Case 2:}} Let $b_{i_t} \in V(\widetilde{F})$ and $b'_t=\widetilde{f}(b_{i_t})$. Then from the proof of Theorem \ref{retraction sequence of blowdown} either $E'_t=\tilde{f}(E_{i_t})$  or $E'_t$ is a face of $\tilde{f}(E_{i_t})$.
Now $\tilde{f}(E_{i_t})$ is either homeomorphic to $E_{i_t}$ as a manifold with corners or homeomorphic to a face of $E_{i_t}$ as a manifold with corners or a blowdown of $E_{i_t}$, see Remark \ref{blowdown of face} and Lemma \ref{lem:blowdown of face}.

{\bf{Subcase 1:}} If $E'_t$ is homeomorphic to $E_{i_t}$ then $$|G_{E'_t}(b'_t)|=|G_{E_{i_t}}(b_{i_t})|.$$

{\bf{Subcase 2:}} If $E'_t $ is homeomorphic to a face of $E_{i_t}$ then, from Proposition \ref{order of group}, $$|G_{E'_t}(b'_t)| ~\text{divides}~ |G_{E_{i_t}}(b_{i_t})|.$$

{\bf{Subcase 3:}} If $E'_t$ is a blowdown of $E_{i_t}$, then from Lemma \ref{projection lemma} and (\ref{equation for blowdown})
\begin{equation}\label{define d_t}
\frac{c_s d_s}{d_{E_{i_t}}} |G_{E'_t}(b'_t)|=|G_{E_{i_t}}(b_{i_t})|~ \text{for some}~ s \in \{1, \dots, (n-d)\}
\end{equation}
where $d_{E_{i_t}}$ comes from (\ref{define d_F}) while computing the determinant of the corresponding matrices given by \ref{singularity 2}. 
Thus $|G_{E'_t}(b'_t)|=\frac{d_{E_{i_t}}}{c_s d_s}|G_{E_{i_t}}(b_{i_t})|$ for some $ s \in \{1, \dots, (n-d)\}$. 
%With $\gcd(|G_{E_i}(b_i)|,p)=1$ for all $i=1, \dots, m$ and $\gcd(\text{denominator of}~ c_s,p)=1$ for $1\leq s \leq (n-d)$
Since we have $d_s \in \ZZ$ from \eqref{Eq_define d_i}, then $d_s$ is a factor of $|G_{E_{i_t}}(b_{i_t})|$. Therefore, if we assume $\gcd(d_{E_{i_t}},p)=1$ then for the above three subcases we have $\gcd(|G_{E'_t}(b'_t)|,p)=1$ where $b'_t=\widetilde{f}(b_{i_t})$.
% which is either in $V(\widetilde{F})$ or adjacent to a vertex of $V(\widetilde{F})$.

{\bf{Case 3:}} Let $b_{i_t}$ be adjacent to a vertex of $V(\widetilde{F})$ and $b'_t=\widetilde{f}(b_{i_t})$ in the blowdown. Then either $E'_t=\tilde{f}(E_{i_t})$  or $E'_t$ is a face of $\tilde{f}(E_{i_t})$. Here also three subcases arise as in \textbf{Case 2} and deduction follows in a similar way. Thus $$\gcd(|G_{E'_t}(b'_t)|,p)=1~\text{for ~all}~1 \leq t \leq m-k.$$ 

The claim $X(Q', \lambda')$ is a blowdown of $X(Q,\lambda)$ follows directly from Proposition \ref{Propn Char vec over Q'}.
\end{proof}

The next two examples show that, in general, we may not relax the hypotheses $(A_2)$ and $(A_3)$ in Theorem \ref{blowdown theorem}.

\begin{example}\label{Eg A2 needed in blowdown}
Let Q be a $3$-dimensional cube and $Q'$ a blowdown of $Q$ as in Figure \ref{Fig_blowdown of a cube}. Define an $\mathcal{R}$-characteristic function $\lambda$ on ${Q}$ by
\begin{align}\label{eq_not_lin_sum}
\lambda(F_0)=(2,1,4), \quad \lambda(F_1)=(6,3,5), \quad
\lambda(F_2)=(3,1,7), \\ \lambda(F_3)=(1,2,6),\quad
\lambda(F_4)=(4,1,3), \quad \lambda(\widetilde{F})=(2,3,5).\nonumber
\end{align}

 Then $|G_Q(b_1)|=5$. Consider the retraction sequence of $Q$  as in Figure \ref{Fig_Figure for necessary condition in blowdown}.

Now we calculate the order of $G_{E_3}(b_3)$. As $E_3$ is the facet $F_1$, we extend $\lambda(F_1)$ to a basis $\{(6,3,5), (1,0,0), (0,2,3)\}$ of $\ZZ^3$. 
Thus the projection map $\rho_{F_1}$ defined in \eqref{eq_proj_map} becomes $$\rho_{F_1} \colon \ZZ^3 \to \ZZ^3/ \left<(6,3,5)\right> \cong \ZZ^2.$$
The facets of $F_1$ are $\{F_1 \cap F_2, F_1 \cap F_3, F_1 \cap F_4, F_1 \cap \widetilde{F}\}$. Therefore the map $\lambda_{F_1} \colon \mathcal{F}(F_1) \to \ZZ^2$ as in \eqref{projection of lambda} is defined by 
\begin{align*}
\lambda_{F_1}(F_1 \cap F_2)&=\rho_{F_1}(\lambda(F_2))=(-63,-16),\\
\lambda_{F_1}(F_1 \cap F_3)&=\rho_{F_1}(\lambda(F_3))=(-35,-8),\\
\lambda_{F_1}(F_1 \cap F_4)&=\rho_{F_1}(\lambda(F_4))=(-14,-4),\\
\lambda_{F_1}(F_1 \cap \widetilde{F})&=\rho_{F_1}(\lambda(\widetilde{F}))=(4,0)
\end{align*}
Thus $|G_{E_3}(b_3)|=64$.
Note that \eqref{eq_not_lin_sum} induces an $\mathcal{R}$-characteristic function on $Q'$ using \eqref{defn of lambda'} though $(2,3,5)$ is not a $\QQ$-linear combination of $(2,1,4)$ and $(6,3,5)$. Then $|G_{Q'}(b'_1)|=7.$
Here, new prime factor 7 arises in the order of singularity at $b_1'$ after blowdown, which was neither in $|G_Q(b_1)|$ nor in $|G_{E_3}(b_3)|$.
Therefore the hypothesis ($A_2$) may not be possible to relax in Theorem \ref{blowdown theorem}.\qed
\end{example}

\begin{example}\label{Eg A3 need in blowdown}
Let $(Q, \lambda)$ and $(Q',\lambda')$ be $\mathcal{R}$-characteristic pairs as in Example \ref{blowdown of cube}. We consider the retraction sequences of $Q$ and $Q'$ as in Figure \ref{Fig_Figure for necessary condition in blowdown}. 
In the induced retraction sequence of $Q'$ from $Q$, $E'_2$ is a blowdown of $E_2$.
Similar calculation to Example \ref{Eg A2 needed in blowdown} gives $|G_{E_2}(b_2)|=1$ but $|G_{E'_2}(b'_2)|=3$.
Here $d_{E_2}=3$ comes while taking determinant, see \eqref{define d_t}.
Thus we cannot relax the hypothesis $(A_3)$ in Theorem \ref{blowdown theorem} in general.\qed
\end{example}

\begin{figure}
\begin{tikzpicture}[scale=.6]
\draw (0,0)--(0,2)--(-1,3)--(-1,1)--cycle;
\draw (0,0)--(2,0)--(2,2)--(0,2);
\draw (2,2)--(1,3)--(-1,3);
\draw [dashed] (-1,1)--(1,1)--(1,3);
\draw [dashed] (1,1)--(2,0);

\node at (2,2) {$\bullet$};
\node [right] at (2,2) {$b_1$};
\node at (1,-.8) {$B_1$};
\draw [->] (3,1.5)--(3.7,1.5);

\begin{scope}[xshift=150]
\draw (0,0)--(0,2)--(-1,3)--(-1,1)--cycle;
\draw [dashed] (-1,1)--(1,1)--(1,3);
\draw [dashed] (1,1)--(2,0);
\draw (1,3)--(-1,3);
\draw (0,0)--(2,0);

\node at (0,2) {$\bullet$};
\node[right] at (0,2) {$b_2$};
\draw [->] (2.5,1.5)--(3.2,1.5);
\node at (1,-.8) {$B_2$};
\end{scope}

\begin{scope}[xshift=300]
\draw (-1,3)--(-1,1)--(0,0);
\draw [dashed] (-1,1)--(1,1)--(1,3);
\draw [dashed] (1,1)--(2,0);
\draw (1,3)--(-1,3);
\draw (0,0)--(2,0);

\node at (2,0) {$\bullet$};
\node[right] at (2,0) {$b_3$};
\draw [->] (2.5,1.5)--(3.2,1.5);
\node at (1,-.8) {$B_3$};
\node [right] at (3.7,1.5) {$\dots$};
\end{scope}

\begin{scope}[yshift=-150]
\draw (0,0)--(-1,3)--(-1,1)--cycle;
\draw (0,0)--(2,0)--(1,3)--(-1,3);
\draw [dashed] (-1,1)--(1,1)--(1,3);
\draw [dashed] (1,1)--(2,0);

\node at (2,0) {$\bullet$};
\node [right] at (2,0) {$b'_1$};
\draw [->] (3,1.5)--(3.7,1.5);
\node at (1,-.8) {$B'_1$};
\end{scope}

\begin{scope}[xshift=150, yshift=-150]
\draw (0,0)--(-1,3)--(-1,1)--cycle;
\draw [dashed] (-1,1)--(1,1)--(1,3);
\draw (1,3)--(-1,3);

\node at (0,0) {$\bullet$};
\node [left] at (0,0) {$b'_2$};
\draw [->] (2.5,1.5)--(3.2,1.5);
\node at (1,-.8) {$B'_2$};
\node [right] at (3.7,1.5) {$\dots$};
\end{scope}

\end{tikzpicture}
\caption{An induced retraction of a blowdown of a cube.}
\label{Fig_Figure for necessary condition in blowdown}
\end{figure}

\begin{remark}
If $c_s$'s are integers in \eqref{lambda after blowdown}, then
$\gcd(\text{denominator of}~ c_s,p)=1$. In this case, we can relax the hypothesis $(A_2)$ in Theorem \ref{blowdown theorem}.
\end{remark}

\begin{theorem}\label{blowdown no torsion}
Let $X(Q', \lambda')$ be a blowdown of $X(Q,\lambda)$ as in \emph{Definition} \ref{blowdown definition} and $(Q,\lambda) $ satisfies the conditions in \emph{Theorem} \ref{blowdown theorem}. Then $H_*(X({Q}',{\lambda}');\mathbb{Z})$ has no $p$-torsion and $H_{\emph{odd}}(X({Q}',{\lambda}');\mathbb{Z}_p)$ is trivial.
\end{theorem}

\begin{proof}
This follows from Theorem \ref{blowdown theorem} and \cite[Theorem 1.1]{BNSS}.
\end{proof}

We note that if $X(Q, \lambda)$ is a quasitoric orbifold, then there is a resolution of singularity $$X(Q(m), \lambda(m)) \to \ldots X(Q(j+1), \lambda(j+1)) \to X(Q(j), \lambda(j)) \to \ldots X(Q(1), \lambda(1))=X(Q,\lambda),$$ where $X(Q(m), \lambda(m))$ is a quasitoric manifold (which is even) and $X(Q(j), \lambda(j)) $ is a blowdown of $ X(Q(j+1), \lambda(j+1))$ for $j=1, \ldots, m$, see \cite[Theorem 2.8]{BSS23}. Then we get the following.

\begin{corollary}\label{blowdown free}  
If $(Q(j),\lambda(j)) $ satisfies the conditions in \emph{Theorem} \ref{blowdown theorem} for any prime $p$ and $j=1, ..., m-1$. Then $H^*(X({Q},{\lambda});\mathbb{Z})$ has no torsion and concentrated in even degrees.
\end{corollary}

\section{The $k$-wedge construction on quasitoric orbifolds and evenness}\label{generalization}

In this section, we introduce $k$-wedge construction on a quasitoric orbifold and show that this gives a new quasitoric orbifold. Next, we show that if the original quasitoric orbifold $X(Q, \lambda)$ satisfies the condition $(A_1)$ as in Theorem \ref{blowdown theorem} then certain $k$-wedge of $X(Q, \lambda)$ satisfies the similar condition. 
Interestingly, this $k$-wedge construction may not be possible to obtained by iterated wedge constructions if $k >1$.

Let $Q$ be an $n$-dimensional simple polytope with $V(Q)=\{v_1, \dots, v_m\}$ and $\mathcal{F}(Q)=\{F^Q_1, \dots, F^Q_r\}$. Then for a $k$-dimensional simplex $\Delta^k$ we get a simple polytope $Q \times \Delta^k$ with $$\mathcal{F}(Q \times \Delta^k)=\{Q \times F^{\Delta}_0, \dots, Q \times F^{\Delta}_k, F^Q_1 \times \Delta^k, \dots, F^Q_r \times \Delta^k\}$$ where $\mathcal{F}(\Delta^k)=\{F_0^{\Delta}, \dots, F_k^{\Delta}\}$. 
Let us consider a blowdown of $Q \times \Delta^k$ of the face $F^Q_{s} \times \Delta$ on $F_{s}^Q$ for some $s \in \{1, \dots, r\}$ and denote it by $(Q \times \Delta^k)^{\prime}$. By Corollary \ref{cor_gen of poly wedge cons}, this is a polytopal $k$-wedge of $Q$ at $F^Q_{s}$. Without loss of generality let $s=r$ and denote the polytopal $k$-wedge by $Q_F(k)$.

Let $\{e_1, \dots, e_k\}$ be the standard basis of $\ZZ^k$. Now we define  a map 
\begin{equation}\label{Eq_char fn on Q Delta}
\widetilde{\lambda} \colon \mathcal{F}(Q \times \Delta^k) \to \ZZ^{n+k},
\end{equation} 
induced from the characteristic function $\lambda$ by the following way
\begin{equation*}
\widetilde{\lambda}(F)=
\begin{cases}
\big(\textbf{0}_k, \lambda(F^Q_j) \big) \quad & \text{ if } F=F^Q_j \times \Delta^k \text{ for } j=1, \dots, r\\
\big(-\sum_{j=1}^k e_j, \lambda(F^Q_{r}) \big) \quad & \text{ if } F= Q \times F^{\Delta}_0\\
\big(1, a, \textbf{0}_{n+k-2}  \big)  \quad & \text{ if } F= Q \times F^{\Delta}_1, ~a \in \ZZ \setminus \{1\}\\
\big(e_j, \textbf{0}_n \big) \quad & \text{ if } F= Q \times F^{\Delta}_j \text{ for } j=2, \dots, k
\end{cases}
\end{equation*}
where $\textbf{0}_j$ represents the zero vector in $j$-dimension, depending on the condition on the facet $F$.

\begin{lemma}\label{lemma on char fn over Q delta}
Let $(Q, \lambda)$ be an $\mathcal{R}$-characteristic pair and $\Delta^k$ a $k$-simplex. Then the map $\widetilde{\lambda}$ defined in \eqref{Eq_char fn on Q Delta} is an $\mathcal{R}$-characteristic map over $Q \times \Delta^k$.
\end{lemma}

\begin{proof}
We investigate the order of singularities defined in \eqref{singularity 2} at the vertices of $Q \times \Delta^k$ and show they are non-zero.
Let $b_i \in V(Q \times {\Delta}^k)$ with $b_i=v_{\ell} \times v^{\Delta}$ for $v_{\ell} \in V(Q)$ and $v^{\Delta} \in V(\Delta^k)$.
If $a=0$ then clearly $\widetilde{\lambda}$ is an $\mathcal{R}$-characteristic function and $|G_{Q \times \Delta^k}(b_i)|=|G_{Q}(v_{\ell})|$.

Now let $a \neq 0,1$. If $b_i \in V(Q \times \Delta^k) \setminus V(Q \times F_1^{\Delta})$ then $$|G_{Q \times \Delta^k}(b_i)|=|G_{Q}(v_{\ell})|$$ where $b_i=v_{\ell} \times v^{\Delta}$. 
Let $b_i \in V(Q \times F_1^{\Delta})$ with $b_i=v_{\ell} \times v^{\Delta}$ and $v_{\ell} = \bigcap_{t=1}^n F_{j_t}^Q$. 
Then
\begin{equation}\label{Eq_bi in product}
b_i=\big(\bigcap_{t=1}^n (F_{j_t}^Q \times \Delta^k)\big) \bigcap \big( \bigcap_{\substack{j=0\\j \neq \alpha}}^k (Q \times F_j^{\Delta}) \big)
\end{equation}
for $\alpha=0,2, \dots, k$. 
To calculate the order of $G_{Q \times \Delta^k}(b_i)$, %in the proof of Lemma \ref{lemma on char fn over Q delta}, 
we can visualise the matrix associated to the vertex $b_i$ in $Q \times \Delta^k$ as the following block matrix $$A_{b_i}^{Q \times \Delta^k}= \begin{pmatrix}
    \textbf{0} & \textbf{B} \\
    A_{v_\ell}^Q & \textbf{C}
\end{pmatrix}$$ where $A_{v_\ell}^Q$ is defined as in \eqref{Eq_associated matrices} and \textbf{B} and \textbf{C} are determined by the vectors assigned to the facets $Q \times F^{\Delta}_j$ of $Q \times \Delta^k$ for $j=0,1, \dots, k$. Thus $$|G_{Q \times \Delta^k}(b_i)|=|\det A_{b_i}^{Q \times \Delta^k}|=|\det A_{v_\ell}^Q |\times| \det \textbf{B}|.$$
If $\alpha=0 \text{ or } 2$ then $|\det \textbf{B}|=1$ and $|G_{Q \times \Delta^k}(b_i)|=|G_{Q}(v_{\ell})| \neq 0$.
If $\alpha \neq 0,2$ then $|\det \textbf{B}|=|(1-a)|$ and $|G_{Q \times \Delta^k}(b_i)|=|(1-a)||G_{Q}(v_{\ell})| \neq 0$.
This concludes the proof of the lemma.

\end{proof}

Note that $|\mathcal{F}(Q_F(k))|=|\mathcal{F}(Q \times \Delta^k)|-1$, since the facet $F_r^Q \times \Delta^k$ is identified with $F^Q_r$ after blowdown.
We recall that the facet set of $Q_F(k)$ is defined in \eqref{Eq_facet set of generalized polytopal wedge}. Now we restrict $\widetilde{\lambda}$ in \eqref{Eq_char fn on Q Delta} to obtain
\begin{equation}\label{Eq_characteristic function on polytopal wedge}
\lambda_F^k \colon \mathcal{F}(Q_F(k)) \to \ZZ^{n+k}
\end{equation}
by the following way
\begin{equation*}
\lambda_F^k(F_i)=
\begin{cases}
\big(\textbf{0}_k, \lambda_i \big) \quad & \text{ for } i=1, \dots, r-1\\
\big(-\sum_{j=1}^k e_j, \lambda(F^Q_{r}) \big) \quad & \text{ for } i=r\\
\big( 1, a, \textbf{0}_{n+k-2} \big) \quad & \text{ for } i=r+1, ~a \in \ZZ \setminus \{1\} \\
\big( e_s, \textbf{0}_n  \big)  \quad & \text{ for } i=r+s \text{ and } s=2, \dots, k
\end{cases}
\end{equation*}
where $\textbf{0}_j$ represents the zero vector of dimension $j$, depending on the condition on the facet $F_i$ of $Q_F(k)$. % and $a \neq 1$. From the construction we have $\lambda_F^k$ is an $\mathcal{R}$-characteristic function over $Q_F(k)$. 

\begin{lemma}
Let $(Q, \lambda)$ be an $\mathcal{R}$-characteristic pair over an $n$-dimensional simple polytope $Q$ with $\mathcal{F}(Q)=\{F_1, \dots, F_r\}$. If $Q_F(k)$ is $k$-wedge of $Q$ at $F=F_r$ and $\lambda_F^k \colon \mathcal{F}(Q_F(k)) \to \ZZ^{n+k}$ is defined as in \eqref{Eq_characteristic function on polytopal wedge}, then $\lambda_F^k$ is an $\mathcal{R}$-characteristic function on $Q_F(k)$.
\end{lemma}

\begin{proof}
This proof is similar to the proof of Lemma \ref{lemma on char fn over Q delta}. 
\end{proof}

\begin{definition}
Let $X(Q,\lambda)$ be a quasitoric orbifold over a simple polytope $Q$ with a facet $F$ and $Q_F(k)$ a polytopal $k$-wedge of $Q$ at $F$. Let $\lambda_F^k$ be defined as in \eqref{Eq_characteristic function on polytopal wedge} induced from $\lambda$.
Then we call the quasitoric orbifold $X(Q_F(k), \lambda_F^k)$ a $k$-wedge of the quasitoric orbifold $X(Q, \lambda)$.
\end{definition}

\begin{remark} 
Observe that $(Q_F(k), \lambda_F^k)$ is a restriction of the characteristic pair $(Q \times \Delta^k , \widetilde{\lambda})$.
Thus the quasitoric orbifold $X(Q_F(k), \lambda_F^k)$ is a blowdown of the quasitoric orbifold $X(Q \times \Delta^k , \widetilde{\lambda})$.
Moreover, if $a=0$ in \eqref{Eq_characteristic function on polytopal wedge} we can use Theorem \ref{blowdown theorem} to the quasitoric orbifold $X(Q_F(k), \lambda_F^k)$ as a blowdown of $X(Q \times \Delta^k , \widetilde{\lambda})$ and get similar conclusion as in Theorem \ref{blowdown no torsion}. Some results for the cases for $a \neq 0$ are discussed further in the following.
\end{remark}

\begin{theorem}\label{proposition on k wedge of quasi orbi}
Let $(Q, \lambda)$ be a quasitoric orbifold over a simple polytope $Q$ with a facet $F$ and for a prime $p$ there exists a retraction sequence $\{(B_{\ell}, E_{\ell}, v_{\ell})\}_{\ell=1}^m$ such that $\{V(F) =\{v_{m-\alpha+1}, \ldots, m\}$ and satisfying $\gcd(|G_{E_{\ell}}(v_{\ell})|,p)=1$ for $\ell=1, \dots, m$.
If $X(Q_F(k), \lambda_F^k)$ is a $k$-wedge of $X(Q, \lambda)$ where $\lambda_F^k$ is defined as in \eqref{Eq_characteristic function on polytopal wedge} such that $\gcd(|1-a|,p)=1$, then there is no $p$-torsion in $H_*(X(Q_F(k), \lambda_F^k);\ZZ)$ and $H_{\text{odd}}(X(Q_F(k), \lambda_F^k); \ZZ_p)=0$.
\end{theorem}

\begin{proof}
Recall the facets of $Q_F(k)$ from \eqref{Eq_facet set of generalized polytopal wedge}, the induced retraction sequence $\{(B'_t, E_t', b_t')\}_{t=1}^u$ of $Q_F(k)$ from Corollary \ref{cor_k wed_ret}, and the $\mathcal{R}$-characteristic vector is defined in \eqref{Eq_characteristic function on polytopal wedge}.
If we prove $\gcd (|G_{E'_{t}}(b'_t)|,p)=1$ for $t=1, \dots, u$, we can conclude the result using \cite[Theorem 1.1]{BNSS}. For that, we have to deal with the following cases. 

{\bf{Case 1:}} Let $1 \leq t \leq (k+1)(m-\alpha)$. In this case $ E'_t = E_{\ell} \times \Delta^{k+1-s}$ and $b'_t = (v_{\ell}, e_{s})$ for $t=(k+1)\ell - (k+1-s)$ where $\ell =1, \dots , m-\alpha$ and $s=1, \dots, k+1$. 
Let $\dim(E_{\ell})=d$. Then $\dim (E'_t)=d+q$ and $$E'_t= (\bigcap_{t=1}^{n-d} (F^Q_{\ell_t} \times \Delta^k)) \bigcap (\bigcap_{j=1}^{k-q} (Q \times F^{\Delta}_{s_j}))$$ where $0 \leq q \leq k$. 
From the discussion in Subsection \ref{Subsec_Basics of quasitoric orbifold}, we obtain a $(d+q) \times (d+q)$ matrix $A^{E'_t}_{b'_t}$ associated to the vertex $b'_t$ in $E'_t$ by projecting $\lambda^k_F$ on the face $E_t'$ as follows.
First we extend the set of $d+q$ vectors $$S({E'_t})=S({E_\ell}) \bigcup S({\Delta^q}):=\{\lambda_F^k(F^Q_{\ell_t} \times \Delta^k) ~|~ t=1, \dots, n-d \} \bigcup \{\lambda_F^k(Q \times F^{\Delta}_{s_j})~|~ j=1, \dots, k-q\}$$ to a basis of $\ZZ^{n+k}$. Since the first $k$ entries of the vectors in $S({E_{\ell}})$ are zeros, we extend them to $n$ linearly independent vectors in $\ZZ^{n+k}$ similar to the extension of $\{\lambda(F^Q_{\ell_t}) ~|~ t=1, \dots, n-d\}$ to a basis of $\ZZ^n$ in $Q$. We denote this linearly independent set of $n$ vectors by $S({Q_\ell})$. 
Also along with $S({\Delta^q})$, we add $q$ many vectors from the standard basis vectors $\{e_1, \dots, e_k\}$ of $\ZZ^{n+k}$ to extend $S({Q_{\ell}})$ to a basis $S({n+k})$ of $\ZZ^{n+k}$.

Now if we visualize the matrix $A^{E'_t}_{b'_t}$ as block matrix of the form $$A^{E'_t}_{b'_t}= \begin{pmatrix}
    \textbf{M}^1_{(q \times d)} & \textbf{M}^2_{(q \times q)} \\
    \textbf{M}^3_{(d \times d)} & \textbf{M}^4_{(d \times q)}
\end{pmatrix}$$ then we have $M^1=\textbf{0}_{(q \times d)}$ and  $M^3=\begin{pmatrix}
    \lambda_{E_{\ell}}(E_{\ell} \cap F_{i_1})^t & \dots & \lambda_{E_{\ell}}(E_{\ell} \cap F_{i_1})^t
\end{pmatrix}_{d \times d}=A^{E_{\ell}}_{v_\ell}$  from the above discussion where $v_\ell= \cap_{j=1}^d (E_{\ell} \cap F_{i_j})$. Thus $$|G_{E'_t}(b'_t)|= |\det A^{E'_t}_{b'_t}|= |\det A^{E_{\ell}}_{v_\ell} | \times |\det M^2|=|G_{E_\ell}(v_\ell)| \times |\det M^2|.$$

If $E'_t \cap F_{r+1} \neq \varnothing$, then there exists two subcases. If $e_2 \in S({n+k})$, then $|\det M^2|=1$. Otherwise, $|\det M^2|=|(1-a)|$. This implies $|G_{E'_t}(b'_t)|$ divides $|(1-a)||G_{E_\ell}(v_\ell)|$.

If $E'_t \cap F_{r+1} = \varnothing$, then $|G_{E'_t}(b'_t)|=|G_{E_\ell}(v_{\ell})|.$ For $\gcd (|1-a|,p)=1$, we can conclude $$\gcd (|G_{E'_{t}}(b'_t)|,p)=1~~~ \text{for} ~t=1, \dots, (k+1)(m-\alpha).$$

{\bf{Case 2:}} Let $(k+1)(m-\alpha)+1 \leq t \leq u=(k+1)(m-\alpha)+\alpha$. In this case $E'_t=E_{m- \alpha +\ell}$ and $b'_t=v_{m- \alpha +\ell}$ for $t=(k+1)(m-\alpha)+\ell$ and $\ell=1, \dots, \alpha$. Then $|G_{E'_t}(b'_t)|=|G_{E_{m-\alpha+\ell}}(v_{m - \alpha +\ell})|$.
Thus, we show $$\gcd (|G_{E'_{t}}(b'_t)|,p)=1~~~ \text{for} ~t=1, \dots, (k+1)(m-\alpha)+\alpha$$ and eventually conclude the result.
\end{proof}

\begin{example}\label{Eg_blowdown can not be obtained by polytopal wedge construction}
In Figure \ref{Fig_Blowdown that can not obtained by polytopal wedge construction.}, we show a blowdown $Q'$ of a simple polytope $Q$ that cannot be obtained by a polytopal wedge construction.
Define an $\mathcal{R}$-characteristic function on $Q$ by
\begin{align}
\lambda(F_0)=(0,2,1), \quad \lambda(F_1)=(1,1,2), \quad \lambda(\widetilde{F})=(1,3,3), \\
\lambda(F_2)=(0,1,1), \quad \lambda(F_3)=(1,0,1), \quad
\lambda(F_4)=(1,0,0), \quad  \lambda(F_5)=(3,2,7). \nonumber
\end{align}
Then $(Q,\lambda)$ is an $\mathcal{R}$-characteristic pair and provides us a quasitoric orbifold $X(Q,\lambda)$.

\begin{figure}
\begin{tikzpicture}[scale=.6]
\draw (0,0)--(0,2)--(-.5,2.5)--(-1,2.5)--(-1,1)--cycle;
\draw[fill=yellow, opacity=.4] (0,0)--(2,0)--(2,2)--(0,2)--cycle;
\draw [thick] (0,0)--(2,0);
\draw (2,2)--(1,3)--(-.5,3)--(-.5,2.5);
\draw (-.5,3)--(-1,2.5);
\draw [dashed] (-1,1)--(1,1)--(1,3);
\draw [dashed] (1,1)--(2,0);

\draw (0,2.5)--(0,3.4);
\node [above] at (0,3.4) {$F_0$};

\draw (.5,0) to [out=240, in=330] (-.5,-.3);
\draw[dotted,thick] (.5,0)--(.7,.3);
\node [left] at (-.5,-.3) {$F_1$};

\draw [dotted, thick] (1.5,1.7)--(2,1.9);
\draw (2,1.9)--(2.5,2.1);
\node [right] at (2.5,2.1) {$F_2$};

\draw[dotted,thick] (-.7,1.8)--(-1,2);
\draw (-1,2)--(-1.7,2.5);
\node [left] at (-1.7,2.5) {$F_3$};

\draw (-.7,1.3)--(-1.7,1.3);
\node [left] at (-1.7,1.3) {$F_4$};

\draw (-.6,2.6)--(-1,3.4);
\node [above] at (-1,3.4) {$F_5$};

\draw (1.7,.7)--(2.3,.7);
\node [right] at (2.3,.7) {$\widetilde{F}$};

\draw [thick] (0,0)--(2,0);
\draw [->] (1.7,-.4)--(1.5,0);
\node [right] at (1.5,-.6) {$F$};

\node at (1,-1) {$Q$};

\draw[->] (4.5,1.5)--(8,1.5);
\node [above] at (6.2,1.5) {$\text{Blowdown}$};
\node [right] at (4.8,1) {$\text{of}~\widetilde{F}~\text{on}~F$};

\begin{scope}[xshift=350]
\draw (0,0)--(-.7,2.2)--(-1,2.5)--(-1,1)--cycle;
\draw (0,0)--(2,0)--(1,3)--(-.5,3)--(-1,2.5);
\draw (-.7,2.2)--(-.5,3);
\draw [dashed] (-1,1)--(1,1)--(1,3);
\draw [dashed] (1,1)--(2,0);

\node at (1,-1) {$Q'$};

\draw (0,2.5)--(0,3.4);
\node [above] at (0,3.4) {$F'_0$};

\draw (.5,0) to [out=240, in=330] (-.5,-.3);
\draw[dotted,thick] (.5,0)--(.7,.3);
\node [left] at (-.5,-.3) {$F'_1$};

\draw [dotted, thick] (1.2,1.3)--(1.5,1.5);
\draw (1.5,1.5)--(2,1.8);
\node [right] at (2,1.8) {$F'_2$};

\draw[dotted,thick] (-.7,1.8)--(-1,2);
\draw (-1,2)--(-1.7,2.5);
\node [left] at (-1.7,2.5) {$F'_3$};

\draw (-.6,1.3)--(-1.2,.5);
\node [left] at (-1.2,.5) {$F'_4$};

\draw (-.7,2.6)--(-1,3.4);
\node [above] at (-1,3.3) {$F'_5$};

\draw [thick] (0,0)--(2,0);
\draw [->] (1.7,-.4)--(1.5,0);
\node [right] at (1.5,-.6) {$F'$};

\end{scope}

\end{tikzpicture}
\caption{A blowdown that cannot be obtained by polytopal $k$-wedge construction.}
\label{Fig_Blowdown that can not obtained by polytopal wedge construction.}
\end{figure} 

We define the $\mathcal{R}$-characteristic function $\lambda'$ on $Q'$ using \eqref{defn of lambda'}.
Then $(Q',\lambda')$ is a restriction of $(Q,\lambda)$. Therefore $X(Q',\lambda')$ is a blowdown of $X(Q,\lambda)$. Note that $X(Q',\lambda')$ cannot be obtained by a $J$-construction of \cite{BSS2} on a $4$-dimensional quasitoric orbifold as its orbit space is a polygon.
\qed
\end{example}

\vspace{1cm}

\noindent {\bf Acknowledgment.}
The authors thank Mainak Poddar and Jongbaek Song for helpful discussions. The first author thanks IIT Madras' for PhD fellowship. The second author thanks `International Office IIT Madras' and `Science and Engineering Research Board India' for research grants. The third author thanks `University Grants Commission of India' for PhD fellowship.

\bibliographystyle{abbrv}
\bibliography{k-wedge}

\end{document}